\documentclass[11pt]{article}
\usepackage[margin=1in]{geometry}
\usepackage[T1]{fontenc}
\usepackage{lmodern}
\usepackage{amsmath,amssymb,amsthm}
\usepackage{graphicx,booktabs}
\usepackage{algorithm,algpseudocode}
\usepackage[numbers,sort&compress]{natbib}
\usepackage[hidelinks]{hyperref}
\newtheorem{theorem}{Theorem}[section]
\newtheorem{lemma}[theorem]{Lemma}
\newtheorem{proposition}[theorem]{Proposition}
\newtheorem{corollary}[theorem]{Corollary}
\newcommand{\R}{\mathbb{R}}

\newcommand{\esssup}{\operatorname*{ess\,sup}}
\newcommand{\essinf}{\operatorname*{ess\,inf}}
\DeclareMathOperator*{\argmin}{arg\,min}
\DeclareMathOperator*{\argmax}{arg\,max}

\newcommand{\FixedFineError}{$0.297\%$}

\newcommand{\ReferenceGridDifference}{$0.00180\%$}
\newcommand{\ReferenceTimeDifference}{$0.00023\%$}

\title{Finite-difference neural policy iteration for high-dimensional degenerate Hamilton--Jacobi equations}
\author{%
  Hee Jun Yang\thanks{National Institute for Mathematical Sciences (NIMS),
    70, Yuseong-daero 1689 beon-gil, Yuseong-gu, Daejeon 34047, Republic of Korea.}
  \and
  Yeoneung Kim\thanks{Department of Industrial Engineering, Yonsei University,
    50 Yonsei-ro, Seodaemun-gu, Seoul 03722, Republic of Korea.
    Corresponding author. Email: \href{mailto:yeoneung@yonsei.ac.kr}{yeoneung@yonsei.ac.kr}.}}
\date{}

\begin{document}
\maketitle

\begin{abstract}
We develop a finite-difference neural policy iteration method for
high-dimensional degenerate Hamilton--Jacobi--Bellman and Isaacs equations.
Artificial viscosity yields a monotone finite-difference operator, and
neural networks approximate the linear equations of semi-discrete policy
iteration using local stencil values. This avoids global tensor-grid storage.
For measurable feedbacks, we bound the propagation of evaluation,
terminal, gradient, and control optimization errors by a geometrically
weighted sum. Discrete comparison, consistency, and vanishing viscosity
relate the approximation to the original degenerate equation.
On coupled manufactured problems with rank-deficient diffusion, the method
achieves lower mean errors than direct physics-informed neural networks
at every tested dimension from five to fifty under equal training time.
At fifty dimensions and $900$ seconds, the mean relative $L^2$ error over
three paired seeds is about $40\%$ lower. Two-dimensional refinement
attains \FixedFineError\ error against the original degenerate solution,
with independent grid references separating approximation error from
viscosity bias. In a twenty-dimensional game with finite control sets,
reusing policies reduces the measured time to a relative error of $5\%$ by
$24.3\%$ compared with direct minimization of the same discrete residual
with fresh samples at every update.
\end{abstract}

\section{Introduction}\label{sec:introduction}

Degenerate Hamilton--Jacobi--Bellman (HJB) and Hamilton--Jacobi--Isaacs
(HJI) equations describe optimal control and differential games when
noise acts in only some state directions or vanishes entirely.
The diffusion operator then lacks a positive ellipticity constant, and
optimal feedbacks may be discontinuous. In differential games, the
Hamiltonian may also be nonconvex in the gradient. Monotone discretizations
converge to the viscosity solution under consistency, stability, and
comparison~\cite{barlessouganidis1991}, but a tensor grid with $M$ points
per coordinate stores $M^d$ values. High-dimensional computation therefore
requires a representation that avoids this storage cost while retaining
a consistent approximation of the degenerate equation.

Neural PDE solvers approximate the solution at sampled states.
The deep Galerkin and deep backward stochastic differential equation
methods established this approach for high-dimensional
problems, including HJB equations
\cite{sirignano2018dgm,han2018}. Physics-informed neural networks (PINNs)
minimize differential-equation residuals evaluated by automatic
differentiation (AD) \cite{raissi2019}.
Numerical differentiation can replace input derivatives in the training
residual \cite{chiu2022can}.
For Hamilton--Jacobi equations, Esteve-Yag\"ue, Tsai, and Massucco
\cite{esteve2025fdleast} evaluate finite-difference stencils with neural
networks and minimize the resulting discrete residual.
Bokanowski, Esteve-Yag\"ue, and Tsai \cite{bokanowski2026residual}
develop error bounds in terms of the residual and monotone approximation theory
for direct neural minimization of the nonlinear discrete residual.
Huang, Kalise, and Zhu~\cite{huang2026residual} use PDE residuals to
select polynomial basis functions for high-dimensional HJB approximation.

Policy iteration (PI) replaces the nonlinear equation by a sequence of
linear policy-evaluation equations and control updates.
Its convergence and sensitivity to inexact evaluations have been
studied for controlled diffusions \cite{kerimkulov2020howard}.
For continuous-time control from discrete data, Zhu~\cite{zhu2026phibe}
develops PDE-based policy evaluation, and Yao et al.~\cite{yao2026phibepi}
quantify how discretization, approximation, and sampling errors
propagate through PI with off-policy data.
Kundu and Kunisch~\cite{kundukunisch2027} analyze bilevel PI for first-
and second-order HJI equations with control constraints and use a
semi-implicit upwind scheme in their numerical tests.
Ito, Reisinger, and Zhang \cite{ito2021neuralpi} analyze inexact neural
PI for stationary uniformly elliptic game equations with measurable
policy selections. Meng et al.\ \cite{meng2024pinnpi} develop
PINN-based PI with convergence and verification. Interior error bounds
for high-dimensional HJB equations are developed in~\cite{kim2026interior}.
Lee and Kim \cite{lee2025operator} combine semi-discrete PI with
DeepONet to learn solution operators across terminal data.
Gomes et al.~\cite{gomes2026newton} analyze Newton and policy iteration
for smooth, possibly nonconvex Hamiltonians and derive residual-based
error estimates for inexact neural Newton corrections.

Our companion paper~\cite{yang2025pinnpi} establishes convergence of
neural PI for nonconvex HJI equations under uniform ellipticity.
That analysis uses parabolic regularity to justify gradient-based
policy improvement. When the physical diffusion is rank deficient or
vanishes, this regularity need not hold, and the value function may lack
a pointwise gradient. The present work addresses this degenerate setting
with a neural approximation suited to high dimensions.

Tang, Tran, and Zhang \cite{tang2025viscosity} make policy improvement
well defined for deterministic HJB equations by using spatial finite
differences and artificial viscosity in a semi-discrete PI scheme.
Guo, Tran, and Zhang \cite{guo2025policy} extend this approach to
nonconvex Hamiltonians with possibly degenerate, control-independent
diffusion. They prove exponential convergence of PI on a monotone
space--time discretization and estimate the discretization error.
Kim et al.~\cite{kim2026monotone} use a related monotone stencil for
neural PI on first-order HJB equations in the whole space. Their fixed-$h$
analysis uses a Markov-chain representation and accounts for errors
from policy evaluation and learned dynamics.

We apply the viscosity-based stencil to nonconvex HJI equations with
possibly degenerate physical diffusion and use a finite-difference (FD)
PINN to approximate each linear policy-evaluation equation, yielding
FD-PINN-PI.
The network supplies values at shifted space--time inputs, and training
minimizes the residual with the policy fixed. The discrete gradient
defines the control update even when the original value function is
nondifferentiable. Local stencil evaluations replace global grid storage;
when control optimization is expensive, selected policies can also be
reused during evaluation.

The analysis accounts for inexact neural evaluation and control updates.
At fixed artificial viscosity $\nu>0$, a geometrically weighted sum of
evaluation, terminal, gradient, and control optimization errors bounds
the squared value and gradient errors. Only measurability of the
feedbacks is required. A discrete comparison estimate relates the value
error to the nonlinear residual, which is controlled by the policy-evaluation
residual and the change in the discrete gradient. Consistency and vanishing
viscosity then connect the approximation to the original equation; for
$C^{1,2}$ solutions, the viscosity bias is $O(\nu)$.

Numerical experiments demonstrate improved accuracy over direct PINN
for degenerate equations in five to fifty dimensions at equal training
time, verify viscosity refinement against two-dimensional grid references,
and quantify the computational benefit of policy reuse for finite control sets.

\section{Finite-difference PINN policy iteration}\label{sec:methods}

We build on the viscosity-based PI of Tang, Tran, and
Zhang~\cite{tang2025viscosity} and its extension by Guo, Tran, and
Zhang~\cite{guo2025policy} to nonconvex, possibly degenerate equations.
Their discrete gradient and artificial viscosity define the policy update
and give a monotone stencil under~\eqref{eq:main_monotonicity}.
We approximate the resulting linear evaluation equations with neural
networks, computing each residual from values at local stencil points.

\subsection{The target equation and its discrete approximation}
We formulate the method on the periodic domain
$\Omega=\mathbb T_\ell^d=(\R/\ell\mathbb Z)^d$, with period $\ell>0$;
the Air3D test uses the bounded domain specified in Section~\ref{sec:air_main}.
Let $A,B$ be compact control spaces and let $\Sigma_j\ge0$ be diagonal
physical variances independent of the controls. For drift $f$, running
cost $c$, terminal data $g$, and terminal time $T>0$, the target is
\begin{equation}\label{eq:main_original}
 v_t+H(t,x,\nabla v)+\frac12\sum_j\Sigma_j\partial_{x_jx_j}v=0,
 \qquad v(T,x)=g(x),
\end{equation}
where
\begin{equation*}
 \begin{aligned}
 L(t,x,p)(a,b)&=f(t,x,a,b)\cdot p+c(t,x,a,b),\\
 H(t,x,p)&=\sup_{b\in B}\inf_{a\in A}L(t,x,p)(a,b).
 \end{aligned}
\end{equation*}
A singleton $B$ gives the minimizing HJB case.

Let $e_j$ be the $j$th coordinate unit vector. For spatial increment $h$
and artificial viscosity $\nu>0$, define
\begin{equation*}
 \begin{aligned}
 D_j^h w&=\frac{w(t,x+he_j)-w(t,x-he_j)}{2h},\\
 \Delta_j^h w&=\frac{w(t,x+he_j)-2w(t,x)+w(t,x-he_j)}{h^2}.
 \end{aligned}
\end{equation*}
The discrete gradient and diffusion operator are
\begin{equation*}
 \nabla_h w=(D_1^h w,\ldots,D_d^h w),\qquad
 D_{\nu,h}w=\frac12\sum_j(\Sigma_j+\nu)\Delta_j^h w.
\end{equation*}
The spatially semi-discrete equation and the corresponding regularized PDE are
\begin{align}
 \partial_t V_h+H(t,x,\nabla_h V_h)+D_{\nu,h}V_h&=0,
 &V_h(T)&=g,\notag\\
 \partial_t v_\nu+H(t,x,\nabla v_\nu)+D_\nu v_\nu&=0,
 &v_\nu(T)&=g,\label{eq:main_target}
\end{align}
with $D_\nu w=\tfrac12\sum_j(\Sigma_j+\nu)\partial_{x_jx_j}w$.
The original data $f,c,\Sigma,g$ remain fixed as $h$ and $\nu$ vary.
We use the spatially semi-discrete equation to define PI and also
discretize time in the neural training residual below.

\subsection{Policy improvement and neural evaluation}
For momentum $p$, the nested feedback is
\begin{equation}\label{eq:main_nested}
 \begin{aligned}
 a_b(p)&\in\argmin_aL(p)(a,b),\\
 b(p)&\in\argmax_bL(p)(a_b(p),b),\qquad a(p)=a_{b(p)}(p).
 \end{aligned}
\end{equation}
The inner minimizer is computed for each candidate $b$ before maximizing
over $b$. The state and time arguments are suppressed. Here momentum denotes
the value-function gradient used in the control update.
Under the assumptions in Section~\ref{sec:main_setting}, the feedback
can be selected measurably.

Let $\widetilde v^n$ denote the neural approximation after policy stage $n$.
Set $p_h^{n-1}=\nabla_h\widetilde v^{n-1}$ and freeze the
feedback $(\alpha^n,\beta^n)=(a(p_h^{n-1}),b(p_h^{n-1}))$.
Policy evaluation approximates the linear equation
\begin{equation*}
 w_t+L(t,x,\nabla_h w)(\alpha^n,\beta^n)+D_{\nu,h}w=0,
 \qquad w(T)=g.
\end{equation*}

The implementation uses the backward time quotient
$D_\tau^-w=[w(t,x)-w(t-\tau,x)]/\tau$ and trains against
\begin{equation}\label{eq:main_fdpi}
 q_n[w]=D_\tau^-w+L(t,x,\nabla_h w)(\alpha^n,\beta^n)+D_{\nu,h}w.
\end{equation}
Training therefore uses space and time differences at sampled points.
For $t<\tau$, the same network
formula supplies values at $t-\tau\in[-\tau,0)$; periodic features
enforce spatial periodicity where applicable.
The trainable representation is $v_\theta=g+(T-t)N_\theta$, where
$N_\theta$ is a neural network with parameters $\theta$. This enforces
the terminal data. The loss is the sample mean of $|q_n[v_\theta]|^2$,
with boundary penalties when required. At the end of stage $n$, its
parameters $\theta_n$ give $\widetilde v^n=v_{\theta_n}$.
Parameter gradients are computed by AD.

\begin{algorithm}[tbp]
\caption{FD-PINN-PI with sampled policy evaluation}\label{alg:integrated}
\begin{algorithmic}[1]
\State Choose $h,\tau,\nu$, a network, a policy interval $K$, and a training budget.
\State Initialize the trainable value network and its frozen copy.
\While{the training budget remains}
 \If{$K$ updates have elapsed since the last policy refresh}
  \State Copy the trainable weights to the frozen network.
 \EndIf
 \State Draw a fresh collocation batch and evaluate $\nabla_h$ of the frozen network.
 \State Compute the nested feedback~\eqref{eq:main_nested} at these momenta.
 \State Update the trainable weights using~\eqref{eq:main_fdpi} and boundary losses.
\EndWhile
\State Evaluate the final network on independent test points.
\end{algorithmic}
\end{algorithm}

The frozen network defines the feedback at each new batch. The index
$n$ counts policy stages; $K$ is the number of optimizer updates per stage.

For expensive control optimization, we retain a collocation batch and
its selected controls for all $K$ updates in a stage, then replace both.
This evaluates~\eqref{eq:main_fdpi} with stored controls and avoids
repeated evaluation of the frozen network. Section~\ref{sec:finite_controls}
uses this implementation for a game with finite control sets.

\subsection{Direct PINN baseline and computational cost}
The primary baseline, direct PINN, minimizes the nonlinear residual
\begin{equation*}
 R_{\rm PINN}[w]=w_t+H(t,x,\nabla w)
                      +\frac12\sum_j\Sigma_j\partial_{x_jx_j}w
\end{equation*}
using input AD. In the refinement and dimension studies, it solves
the original PDE, while FD-PINN-PI uses the regularization in~\eqref{eq:main_fdpi}.
Both methods are tested against the same original solution and use paired
initializations, the same architecture, precision, sampling distribution,
and optimizer. The Air3D comparison uses the same positive viscosity
for both methods.

Additional comparisons distinguish discretization from policy iteration.
FD-direct minimizes the nonlinear FD residual with the same
stencil, while AD-PINN-PI evaluates frozen policies using input AD.
The policy-interval and finite-control comparisons use FD-direct;
Section~\ref{sec:air_main} compares all four methods for Air3D.

For an interior batch of size $N_b$, the trainable FD residual uses
$(2d+2)N_b$ inputs; Algorithm~\ref{alg:integrated} uses another $2dN_b$
for the frozen FD momentum. Retaining the batch and selected controls
avoids these additional queries and requires one control search per
$K$ updates, compared with one per update for FD-direct.
The number of local queries per sampled state thus grows linearly
with dimension. We measure cost by training time and peak GPU allocation.
Direct PINN computes the input gradient and the diagonal second
derivatives in the active diffusion directions. The implementation
omits zero-coefficient second derivatives.

Training time includes sampling, policy refreshes, residual evaluation,
backpropagation, and optimizer updates. Diagnostics and checkpoint I/O
are excluded. The fixed-target and dimension studies use an RTX~5070;
peak PyTorch CUDA allocation is measured during training. Final
checkpoints follow the prescribed update or elapsed-time budget.
\section{Error estimates}\label{sec:stability}

We first bound the propagation of policy-evaluation errors at fixed
$\nu>0$. Consistency relates the discrete training residual to these
continuum errors. Discrete comparison and vanishing viscosity then give
error estimates relative to the original degenerate solution $v$.

\subsection{Assumptions and approximation errors}\label{sec:main_setting}
The frozen drift and running cost must be measurable for the variational
formulation below. The following result shows how to choose the nested
controls in~\eqref{eq:main_nested} so that this holds for any measurable
momentum field, even when the optimizers are nonunique.

\begin{proposition}\label{prop:main_selector}
Let $A,B$ be nonempty compact metric spaces. If $L$ is Borel measurable
in $(t,x,p)$ and continuous in $(a,b)$, the selectors in
\eqref{eq:main_nested} can be chosen measurably. Their composition with
any measurable momentum field satisfies $L(p)(a(p),b(p))=H(p)$.
\end{proposition}

\begin{proof}
The measurable maximum theorem gives a measurable inner minimizer.
The marginal $m(t,x,p,b)=\min_a L(t,x,p)(a,b)$ is measurable in
$(t,x,p)$ and continuous in $b$ by Berge's theorem~\cite{berge1963}. Applying the
measurable maximum theorem to $m$ on $B$ gives a measurable outer
maximizer~\cite[Theorem 18.19]{aliprantisborder2006}. Composing these
two selections gives a measurable feedback, and $L(p)(a(p),b(p))=H(p)$
follows from the two optimizations.
\end{proof}

The energy analysis concerns the regularized equation at $\nu>0$;
the original physical diffusion may be degenerate.
Let $\Omega=\mathbb T_\ell^d$, $Q_T=(0,T)\times\Omega$, and write $\|\cdot\|$
for the spatial $L^2(\Omega)$ norm. Assume that $f,c$ are bounded Borel
functions, continuous in the compact controls, $g$ is bounded and Borel,
and $\Sigma_j\in W^{1,\infty}(Q_T)$. Set
\begin{equation*}
 \begin{split}
 F&=\esssup_{Q_T}\sup_{a,b}|f(t,x,a,b)|,\qquad
 K_\Sigma=\esssup_{Q_T}\Big(\sum_j|\partial_{x_j}\Sigma_j|^2\Big)^{1/2},\\
 \lambda&=\min_j\essinf_{Q_T}(\Sigma_j+\nu)>0.
 \end{split}
\end{equation*}
Set $\Lambda=\max_j\esssup_{Q_T}(\Sigma_j+\nu)$ and use the energy space
\[
 W(0,T)=\{w\in L^2(0,T;H^1(\Omega)):
 w_t\in L^2(0,T;H^{-1}(\Omega))\}.
\]
This space embeds continuously in $C([0,T];L^2(\Omega))$.
The nondivergence diffusion is interpreted variationally through
\begin{equation}\label{eq:main_divergence}
 D_\nu w=\frac12\sum_j\partial_{x_j}((\Sigma_j+\nu)\partial_{x_j}w)
         -\frac12\sum_j(\partial_{x_j}\Sigma_j)\partial_{x_j}w.
\end{equation}
\begin{proposition}\label{prop:main_wellposed}
Under these coefficient assumptions and for $g\in L^2(\Omega)$, each
measurable frozen-policy equation
\begin{equation}\label{eq:main_wellposed_frozen}
 w_t+f_{\alpha,\beta}\cdot\nabla w+D_\nu w+c_{\alpha,\beta}=0,
 \qquad w(T)=g,
\end{equation}
has a unique solution in $W(0,T)$, as does the modified nonlinear
equation~\eqref{eq:main_target}. The frozen equation also admits any
$L^2(Q_T)$ forcing. Here $f_{\alpha,\beta}$ and $c_{\alpha,\beta}$ denote
the coefficients evaluated at the measurable feedback pair.
\end{proposition}

The proof in \ref{app:wellposed} constructs the linear solution
by coercive forms and the nonlinear solution by a weighted contraction.

Take $\widetilde v^0\in L^2(0,T;H^1(\Omega))$ and
$\widetilde v^n\in W(0,T)$ for $n\ge1$. For measurable policy momenta
$p^{n-1}\in L^2(Q_T;\mathbb R^d)$ and measurable feedbacks, define
\begin{align*}
 r_n&=\partial_t\widetilde v^n+
   L(t,x,\nabla\widetilde v^n)(\alpha^n,\beta^n)+D_\nu\widetilde v^n,
   &\text{evaluation defect},\\
 e_T^n&=\widetilde v^n(T)-g,
   &\text{terminal defect},\\
 \eta_{n-1}&=p^{n-1}-\nabla\widetilde v^{n-1},
   &\text{momentum defect}.
\end{align*}
The control optimization error is bounded by a measurable nonnegative function satisfying
\begin{equation}\label{eq:main_attainment}
 |L(t,x,p^{n-1})(\alpha^n,\beta^n)-H(t,x,p^{n-1})|
       \le e_{\rm opt}^n(t,x).
\end{equation}
Require $r_n,\eta_{n-1},e_{\rm opt}^n$ to belong to the corresponding
space--time $L^2$ spaces; $e_T^n$ is the $L^2$ trace. In exact nested
improvement $e_{\rm opt}^n=0$. For AD momentum $\eta_{n-1}=0$ in exact
arithmetic, whereas FD momentum contributes its stencil discrepancy.

The approximation therefore solves the perturbed equation
\begin{equation}\label{eq:main_perturbed}
 \partial_t\widetilde v^n+L(t,x,\nabla\widetilde v^n)(\alpha^n,\beta^n)
 +D_\nu\widetilde v^n=r_n,\qquad \widetilde v^n(T)=g+e_T^n.
\end{equation}

\subsection{Linear stability and policy mismatch}
We first estimate the effect of perturbations in one policy evaluation.
The policy-mismatch bound then relates successive evaluations.
\begin{lemma}\label{lem:analysis_energy}
Let $(\alpha,\beta)$ be measurable feedbacks and let $w\in W(0,T)$ solve
\begin{equation}\label{eq:linear_error_equation}
 \partial_tw+f_{\alpha,\beta}\cdot\nabla w
 +\frac12\sum_{j=1}^d(\Sigma_j+\nu)\partial_{x_jx_j}^2w
 =F_0
 \quad\text{in }Q_T,
 \qquad
 w(T)=w_T,
\end{equation}
where $F_0\in L^2(Q_T)$ and $w_T\in L^2(\Omega)$. Define
\[
 \kappa_{\mathrm{fp}}
 :=1+\frac{(2F+K_\Sigma)^2}{\lambda},
 \qquad
 \omega_{\mathrm{fp}}(t):=e^{-\kappa_{\mathrm{fp}}(T-t)}.
\]
Then
\begin{equation}\label{eq:frozen_energy_estimate}
 \omega_{\mathrm{fp}}(t)\|w(t)\|^2
 +\frac{3\lambda}{4}\int_t^T\omega_{\mathrm{fp}}(s)\|\nabla w(s)\|^2\,ds
 \le
 \|w_T\|^2+\int_t^T\omega_{\mathrm{fp}}(s)\|F_0(s)\|^2\,ds.
\end{equation}
\end{lemma}

\begin{proof}
For $w\in W(0,T)$, the map $t\mapsto\|w(t)\|^2$ is absolutely continuous with derivative $2\langle\partial_tw(t),w(t)\rangle$ \cite{lionsmagenes1972}. Testing~\eqref{eq:linear_error_equation} by $w$ and integrating by parts on the torus, using~\eqref{eq:main_divergence}, yields
\[
 -\frac{d}{dt}\|w\|^2
 +\lambda\|\nabla w\|^2
 \le
 2\|F_0\|\|w\|
 +(2F+K_\Sigma)\|\nabla w\|\|w\|.
\]
The constant $K_\Sigma$ arises from Cauchy--Schwarz applied to the coefficient-gradient product. In particular, the term
\[\sum_j(\partial_{x_j}\Sigma_j)(\partial_{x_j}w)w\] introduces no additional dimension factor. Young's inequality gives
\[
 (2F+K_\Sigma)\|\nabla w\|\|w\|
 \le\frac{\lambda}{4}\|\nabla w\|^2
 +\frac{(2F+K_\Sigma)^2}{\lambda}\|w\|^2,
\]
\[
 2\|F_0\|\|w\|\le\|F_0\|^2+\|w\|^2.
\]
Multiplication by $\omega_{\mathrm{fp}}$ and integration from $t$ to $T$ prove~\eqref{eq:frozen_energy_estimate}.
\end{proof}
The Lipschitz estimates below adapt the policy-comparison argument of
Guo, Tran, and Zhang~\cite[Lemmas 2.4--2.5]{guo2025policy} to updates
satisfying~\eqref{eq:main_attainment}. They bound the error through Hamiltonian values
without differentiating the feedback selectors.
\begin{lemma}\label{lem:analysis_H_lip}
For a.e.\ $(t,x)$ and all $p,q\in\R^d$,
\begin{equation}\label{eq:L_H_lipschitz}
\begin{aligned}
\sup_{a,b}|L(t,x,p)(a,b)-L(t,x,q)(a,b)|&\le F|p-q|,\\
|H(t,x,p)-H(t,x,q)|&\le F|p-q|.
\end{aligned}
\end{equation}
\end{lemma}

\begin{proof}
The first inequality follows from the affine dependence $f\cdot p$. The operations $\inf$ and $\sup$ are nonexpansive under uniform perturbations, which gives the second inequality.
\end{proof}

\begin{lemma}\label{lem:analysis_mismatch}
Let $p:Q_T\to\R^d$ be measurable and let $(\alpha,\beta)$ satisfy
\[
 |L(t,x,p)(\alpha,\beta)-H(t,x,p)|\le e_{\rm opt}
\]
a.e. Then for any measurable momentum fields $q,r$,
\begin{equation}\label{eq:nested_policy_mismatch}
 |L(t,x,q)(\alpha,\beta)-H(t,x,r)|
 \le F|q-p|+F|p-r|+e_{\rm opt}
 \quad\text{a.e.}
\end{equation}
\end{lemma}

\begin{proof}
Insert $L(t,x,p)(\alpha,\beta)$ and $H(t,x,p)$ between the two terms and use~\eqref{eq:L_H_lipschitz} together with the assumed optimization defect.
\end{proof}

Set
\[
 u_n:=\widetilde v^n-v_{\nu}.
\]
Applying \eqref{eq:nested_policy_mismatch} with
$p=p^{n-1}$, $q=\nabla\widetilde v^n$, and $r=\nabla v_{\nu}$ gives
\begin{equation}\label{eq:main_mismatch}
\begin{aligned}
&\left|L(t,x,\nabla\widetilde v^n)(\alpha^n,\beta^n)-H(t,x,\nabla v_{\nu})\right|\\
&\quad\le F|\nabla u_n|+2F|\nabla u_{n-1}|+2F|\eta_{n-1}|+e_{\rm opt}^n.
\end{aligned}
\end{equation}

\subsection{Accumulation of policy-iteration errors}
We now combine these estimates to track the contribution of each
policy stage to the error in later iterates.
For $0<\rho<1$, let
\begin{equation*}
 \kappa_\rho=3+\frac{(2F+K_\Sigma)^2}{\lambda}
                    +\frac{16F^2}{3\rho\lambda},\qquad
  \omega_\rho(t)=e^{-\kappa_\rho(T-t)}.
\end{equation*}
Set
\begin{equation*}
 \begin{aligned}
 S_n={}&\|e_T^n\|^2+\|r_n\|_{L^2(Q_T)}^2\\
       &+4F^2\|\eta_{n-1}\|_{L^2(Q_T)}^2
          +\|e_{\rm opt}^n\|_{L^2(Q_T)}^2.
 \end{aligned}
\end{equation*}
For the error $u_n$ defined above, set
$\Phi_0=\lambda\int_0^T\omega_\rho\|\nabla u_0\|^2\,dt$.

\begin{theorem}[Error propagation in inexact policy iteration]\label{thm:main_history}
Under these assumptions, for every $n\ge1$,
\begin{equation}\label{eq:main_history}
 \underbrace{\sup_{t\in[0,T]}\|u_n(t)\|^2+
       \lambda\int_0^T\|\nabla u_n(t)\|^2\,dt}_{E_n}
 \le \frac73 e^{\kappa_\rho T}
       \left(\sum_{m=1}^n\rho^{n-m}S_m+\rho^n\Phi_0\right).
\end{equation}
The estimate holds for measurable, possibly nonunique feedback selectors.
The initialization term is bounded in terms of the data by
\begin{equation}\label{eq:main_initial_bound}
 \Phi_0\le 2\lambda\int_0^T\omega_\rho
                 \|\nabla\widetilde v^0\|^2\,dt
       +\frac83\bigl(\|g\|^2+T|\Omega|C_c^2\bigr),
\end{equation}
where $C_c=\esssup_{Q_T}\sup_{a,b}|c(t,x,a,b)|$.
\end{theorem}

\begin{proof}
Subtract~\eqref{eq:main_target} from~\eqref{eq:main_perturbed}. With
\[
 \ell_n
 :=L(t,x,\nabla\widetilde v^n)(\alpha^n,\beta^n)
 -H(t,x,\nabla v_{\nu}),
\]
the error satisfies
\[
 \partial_tu_n
 +\frac12\sum_{j=1}^d(\Sigma_j+\nu)\partial_{x_jx_j}^2u_n
 =r_n-\ell_n,
 \qquad
 u_n(T)=e_T^n.
\]
Testing by $u_n$ and integrating by parts gives
\begin{align*}
 -\frac{d}{dt}\|u_n\|^2+\lambda\|\nabla u_n\|^2
 &\le
 2\|r_n\|\|u_n\|
 +2\int_\Omega|\ell_n||u_n|\,dx
 +K_\Sigma\|\nabla u_n\|\|u_n\|.
\end{align*}
Using~\eqref{eq:main_mismatch},
\begin{align*}
 2\int_\Omega|\ell_n||u_n|\,dx
 &\le
 2F\|\nabla u_n\|\|u_n\|
 +4F\|\nabla u_{n-1}\|\|u_n\|\\
 &\quad
 +4F\|\eta_{n-1}\|\|u_n\|
 +2\|e_{\rm opt}^n\|\|u_n\|.
\end{align*}
Young's inequality yields
\[
 (2F+K_\Sigma)\|\nabla u_n\|\|u_n\|
 \le
 \frac{\lambda}{4}\|\nabla u_n\|^2
 +\frac{(2F+K_\Sigma)^2}{\lambda}\|u_n\|^2,
\]
\[
 4F\|\nabla u_{n-1}\|\|u_n\|
 \le
 \frac{3\rho\lambda}{4}\|\nabla u_{n-1}\|^2
 +\frac{16F^2}{3\rho\lambda}\|u_n\|^2,
\]
\[
 2\|r_n\|\|u_n\|\le\|r_n\|^2+\|u_n\|^2,
\]
\[
 4F\|\eta_{n-1}\|\|u_n\|
 \le4F^2\|\eta_{n-1}\|^2+\|u_n\|^2,
\qquad
 2\|e_{\rm opt}^n\|\|u_n\|
 \le\|e_{\rm opt}^n\|^2+\|u_n\|^2.
\]
Therefore
\begin{align*}
 -\frac{d}{dt}\|u_n\|^2
 +\frac{3\lambda}{4}\|\nabla u_n\|^2
 &\le
 \kappa_\rho\|u_n\|^2
 +\frac{3\rho\lambda}{4}\|\nabla u_{n-1}\|^2
 \\
 &\quad
 +\|r_n\|^2+4F^2\|\eta_{n-1}\|^2
 +\|e_{\rm opt}^n\|^2.
\end{align*}
Multiply by $\omega_\rho$, integrate from $t$ to $T$, and use the terminal mismatch. We obtain
\begin{align}
 \omega_\rho(t)\|u_n(t)\|^2
 +\frac{3\lambda}{4}\int_t^T\omega_\rho\|\nabla u_n\|^2\,ds
 \le
 S_n+\frac{3\rho\lambda}{4}
 \int_0^T\omega_\rho\|\nabla u_{n-1}\|^2\,ds.
 \label{eq:weighted_geometric_recursion}
\end{align}
Set
\[
 \Phi_n:=\lambda\int_0^T\omega_\rho\|\nabla u_n\|^2\,dt.
\]
Taking $t=0$ in~\eqref{eq:weighted_geometric_recursion} and discarding the nonnegative value term gives
\[
 \Phi_n\le\frac{4}{3}S_n+\rho\Phi_{n-1}.
\]
Hence
\begin{equation}\label{eq:phi_unrolled}
 \Phi_n\le\frac{4}{3}\sum_{m=1}^n\rho^{n-m}S_m+\rho^n\Phi_0.
\end{equation}
Taking the supremum in~\eqref{eq:weighted_geometric_recursion} and inserting~\eqref{eq:phi_unrolled} at level $n-1$ gives
\[
 \sup_t\omega_\rho(t)\|u_n(t)\|^2
 \le
 \sum_{m=1}^n\rho^{n-m}S_m+\frac{3}{4}\rho^n\Phi_0.
\]
Adding this estimate to~\eqref{eq:phi_unrolled}, using $\omega_\rho\ge e^{-\kappa_\rho T}$, and bounding both coefficients by $7/3$ proves~\eqref{eq:main_history}.
To verify the initialization bound~\eqref{eq:main_initial_bound}, test the
modified equation using $|H(t,x,p)|\le F|p|+C_c$. The energy calculation
used for~\eqref{eq:frozen_energy_estimate} gives
\[
 \lambda\int_0^T\omega_{\rm fp}\|\nabla v_\nu\|^2\,dt
 \le\frac43\bigl(\|g\|^2+T|\Omega|C_c^2\bigr).
\]
Since $\omega_\rho\le\omega_{\rm fp}$, applying
$|a-b|^2\le2|a|^2+2|b|^2$ proves the stated bound.
\end{proof}

For fixed $\nu$ and $\rho$, $S_m\le S$ implies
$\limsup_n E_n\le\tfrac73e^{\kappa_\rho T}S/(1-\rho)$,
and $S_n\to0$ implies $E_n\to0$. For the latter, split the sum into
finitely many terms and a tail bounded by $\sup_{m\ge M}S_m/(1-\rho)$.
A perturbation at stage $m$ contributes $\rho^{n-m}S_m$ to the sum.
The bound in~\cite{yang2025pinnpi} uses the maximum residual over previous
iterations; \eqref{eq:main_history} records each stage separately and
allows measurable, nonunique feedbacks.
The prefactor $e^{\kappa_\rho T}$ grows exponentially as $\lambda\to0$,
so Theorem~\ref{thm:main_history} is a fixed-viscosity stability estimate.
Proposition~\ref{prop:main_discrete} instead bounds the discrete value error
by the accumulated global sup-norm residual, without an inverse-ellipticity
prefactor. Section~\ref{sec:main_defects} tests
\eqref{eq:main_history} at a moderate viscosity.

\subsection{Finite-difference consistency}
The network is trained on a discrete residual, whereas the preceding
estimate uses continuum defects. We therefore relate the two by consistency.
Write $\mathcal R_\nu^{\alpha,\beta}[w]
=w_t+L(t,x,\nabla w)(\alpha,\beta)+D_\nu w$ and let
$\mathcal R_{\tau,h,\nu}^{\alpha,\beta}$ denote its time-quotient and
spatial finite-difference approximation in~\eqref{eq:main_fdpi}.
\begin{lemma}\label{lem:analysis_cd}
Let $w\in C^2([-\tau,T];C^4(\Omega))$. Then, for every $j$,
\begin{align}
 |D_j^hw-\partial_{x_j}w|
 &\le\frac{h^2}{6}\|\partial_{x_j}^3w\|_{L^\infty(Q_T)},
 \label{eq:first_difference_consistency}\\
 |\Delta_j^hw-\partial_{x_jx_j}^2w|
 &\le\frac{h^2}{12}\|\partial_{x_j}^4w\|_{L^\infty(Q_T)},\notag\\
 |D_{\tau}^-w-\partial_tw|
 &\le\frac{\tau}{2}\|\partial_{tt}w\|_{L^\infty((-\tau,T)\times\Omega)}.
 \label{eq:time_difference_consistency}
\end{align}
Consequently, for every measurable feedback pair $(\alpha,\beta)$,
the following bound holds almost everywhere in $Q_T$:
\begin{equation}\label{eq:operator_consistency_exact}
 \left|
 \mathcal R_{\nu}^{\alpha,\beta}[w]
 -\mathcal R_{\tau,h,\nu}^{\alpha,\beta}[w]
 \right|
 \le \mathfrak c_{\tau,h}[w],
\end{equation}
where
\begin{align*}
 \mathfrak c_{\tau,h}[w]
 &:=
 \frac{\tau}{2}\|\partial_{tt}w\|_{L^\infty((-\tau,T)\times\Omega)}
 +\frac{Fh^2}{6}\sum_{j=1}^d\|\partial_{x_j}^3w\|_{L^\infty(Q_T)}
 \\
 &\quad
 +\frac{\Lambda h^2}{24}\sum_{j=1}^d\|\partial_{x_j}^4w\|_{L^\infty(Q_T)}.
\end{align*}
\end{lemma}

\begin{proof}
Taylor's formula with integral remainder gives~\eqref{eq:first_difference_consistency}--\eqref{eq:time_difference_consistency}. Moreover,
\[
 |L(t,x,\nabla w)(\alpha,\beta)-L(t,x,\nabla_hw)(\alpha,\beta)|
 \le F|\nabla w-\nabla_hw|,
\]
and
\[
 \frac12\sum_{j=1}^d(\Sigma_j+\nu)
 |\partial_{x_jx_j}^2w-\Delta_j^hw|
 \le\frac{\Lambda h^2}{24}
 \sum_{j=1}^d\|\partial_{x_j}^4w\|_{L^\infty(Q_T)}.
\]
Adding the time, first-derivative, and second-derivative defects proves~\eqref{eq:operator_consistency_exact}.
\end{proof}

Apply these estimates to a neural approximation
$\widetilde v^n\in C^2([-\tau,T];C^4(\Omega))$ with FD momentum
$p^n=\nabla_h\widetilde v^n$. For the momentum defect $\eta_n$ defined in
Section~\ref{sec:main_setting}, \eqref{eq:first_difference_consistency} gives
\begin{equation}\label{eq:gradient_stencil_bound}
 \|\eta_n\|_{L^2(Q_T)}
 \le
 \frac{h^2|Q_T|^{1/2}}{6}
 \left(\sum_{j=1}^d
 \|\partial_{x_j}^3\widetilde v^n\|_{L^\infty(Q_T)}^2
 \right)^{1/2}.
\end{equation}
Let
\[
 q_n:=\mathcal R_{\tau,h,\nu}^{\alpha^n,\beta^n}[\widetilde v^n].
\]
Then~\eqref{eq:operator_consistency_exact} gives
\begin{equation}\label{eq:residual_consistency_l2}
 \|r_n\|_{L^2(Q_T)}
 \le
 \|q_n\|_{L^2(Q_T)}
 +|Q_T|^{1/2}\mathfrak c_{\tau,h}[\widetilde v^n].
\end{equation}
Independent validation with a global residual bound gives
the sampling estimate in Lemma~\ref{lem:analysis_sampling}.
Equations~\eqref{eq:gradient_stencil_bound} and~\eqref{eq:residual_consistency_l2}
bound the momentum and evaluation terms in $S_n$.

\subsection{Discrete residual estimates}\label{sec:main_discrete_bridge}

A second estimate applies directly to the difference scheme.
Monotonicity and invariance under addition of constants yield a discrete
comparison estimate, as in~\cite{guo2025policy,lee2025operator,bokanowski2026residual}.

Let $\Omega=\mathbb T_\ell^d$, $h>0$, $t_k=k\tau$, and $N\tau=T$.
Here $k$ indexes time levels; the policy stage is fixed.
Suppose $f,c,\Sigma_j$ are bounded, periodic, continuous functions of their
state/time arguments, jointly continuous in the compact controls where
applicable, and the physical variances satisfy $\Sigma_j\ge0$.
Diffusion is diagonal and independent of controls. Let $g$ be bounded and
periodic, and write $a_j^k(x)=\Sigma_j(t_k,x)+\nu$. Assume
\begin{equation}\label{eq:main_monotonicity}
 a_j^k(x)\ge h|f_j(t_k,x,a,b)|\quad\hbox{for every }x,a,b,j,k,
 \qquad\frac{\tau}{h^2}\sum_j a_j^k(x)\le1.
\end{equation}
All sup norms in this section are taken over $\Omega$.
The same argument applies on a periodic grid closed under
the stencil shifts, with the norm then restricted to that grid.

Write $D_{\nu,h}^kz=\tfrac12\sum_j a_j^k\Delta_j^hz$ for the diffusion
operator at $t_k$, and define
$\mathcal M_k[z]=z+\tau H(t_k,x,\nabla_h z)+\tau D_{\nu,h}^k z$.
For a fixed control pair, its affine stencil is
\begin{equation}\label{eq:supp_discrete_weights}
 \begin{split}
 \mathcal M_k^{a,b}[z](x)
 &=\left(1-\frac{\tau}{h^2}\sum_j a_j^k(x)\right)z(x)+\tau c(t_k,x,a,b)\\
 &\quad+\sum_j\frac{\tau}{2h^2}
       \left(a_j^k(x)+h f_j(t_k,x,a,b)\right)z(x+he_j)\\
 &\quad+\sum_j\frac{\tau}{2h^2}
       \left(a_j^k(x)-h f_j(t_k,x,a,b)\right)z(x-he_j).
 \end{split}
\end{equation}
Because the center and diffusion terms are independent of the controls,
$\mathcal M_k[z]=\sup_b\inf_a\mathcal M_k^{a,b}[z]$.
Every coefficient of $z$ in~\eqref{eq:supp_discrete_weights} is nonnegative
under~\eqref{eq:main_monotonicity}, and their sum is one.
Thus $z\le y$ implies $\mathcal M_k[z]\le\mathcal M_k[y]$ and
$\mathcal M_k[z+C]=\mathcal M_k[z]+C$ for spatial constants $C$.
In particular,
\begin{equation}\label{eq:supp_nonexpansive}
 \|\mathcal M_k[z]-\mathcal M_k[y]\|_\infty\le\|z-y\|_\infty.
\end{equation}
Define the exact discrete solution by
$V_{h,\tau}^N=g$, $V_{h,\tau}^{k-1}=\mathcal M_k[V_{h,\tau}^k]$.
In particular, $\|V_{h,\tau}^k\|_\infty\le\|g\|_\infty+(T-t_k)C_c$,
where $C_c=\sup_{t,x,a,b}|c(t,x,a,b)|$.

\begin{proposition}[Discrete residual error bound]\label{prop:main_discrete}
For any sequence of bounded periodic functions $w^0,\ldots,w^N$, set
\[
 d^k=\frac{w^k-w^{k-1}}{\tau}+H(t_k,x,\nabla_h w^k)
                       +D_{\nu,h}^k w^k,\qquad 1\le k\le N.
\]
Under~\eqref{eq:main_monotonicity},
\begin{equation}\label{eq:main_discrete_certificate}
 \|w^k-V_{h,\tau}^k\|_\infty
 \le\|w^N-g\|_\infty+\tau\sum_{j=k+1}^N\|d^j\|_\infty.
\end{equation}
For an inexact policy evaluation, suppose the frozen policy satisfies
\[
 |L(t_k,x,p^{k})(\alpha^k,\beta^k)-H(t_k,x,p^{k})|\le\epsilon^k(x),
\]
and write $q^k$ for the same time/spatial residual with
that policy frozen. If $F=\sup_{t,x,a,b}|f(t,x,a,b)|$, then
\begin{equation}\label{eq:main_discrete_mismatch}
 \|d^k\|_\infty\le\|q^k\|_\infty
              +2F\|\nabla_h w^k-p^k\|_\infty+\|\epsilon^k\|_\infty.
\end{equation}
\end{proposition}

\begin{proof}
Rearranging the definition gives
$w^{k-1}=\mathcal M_k[w^k]-\tau d^k$.
Subtracting the exact recursion and using~\eqref{eq:supp_nonexpansive} yields
\[
 \|w^{k-1}-V_{h,\tau}^{k-1}\|_\infty
 \le\|w^k-V_{h,\tau}^k\|_\infty+\tau\|d^k\|_\infty.
\]
Backward induction proves~\eqref{eq:main_discrete_certificate}.
For the policy estimate, the time and diffusion terms cancel in $d^k-q^k$.
Insert $H(t_k,x,p^k)$ and $L(t_k,x,p^k)(\alpha^k,\beta^k)$.
The two Lipschitz estimates of Lemma~\ref{lem:analysis_H_lip}, with the
constant $F$ defined above, together with the control optimization
bound, give
$|d^k-q^k|\le2F|\nabla_h w^k-p^k|+\epsilon^k$ pointwise.
Taking sup norms proves~\eqref{eq:main_discrete_mismatch}.
\end{proof}

The stencils in Sections~\ref{sec:dimension_scaling},
\ref{sec:fixed_target}, and~\ref{sec:finite_controls}
satisfy~\eqref{eq:main_monotonicity}. The experiments report sampled residual norms.

\subsection{Accuracy for the original degenerate equation}\label{sec:main_total_error}
It remains to control the difference between the discrete target and
the solution of the original degenerate equation.
At $t_k$, the triangle inequality and~\eqref{eq:main_discrete_certificate} give
\begin{equation}\label{eq:main_total_discrete}
\begin{aligned}
 \|w^k-v(t_k)\|_\infty
 &\le\|w^N-g\|_\infty+\tau\sum_{j=k+1}^N\|d^j\|_\infty\\
 &\quad+\|V_{h,\tau}^k-v(t_k)\|_\infty.
\end{aligned}
\end{equation}
For the last term, assume that $g\in C(\Omega)$, the fixed coefficients
satisfy the continuity hypotheses of Section~\ref{sec:main_discrete_bridge},
and the original terminal-value problem has a continuous viscosity solution
$v$ and comparison for bounded semicontinuous sub- and supersolutions.
Comparison is understood on $[0,T]\times\Omega$, with one-sided time
tests at $t=0$.
Consider $h,\tau,\nu\to0$ with $0<\nu\le1$, $T/\tau\in\mathbb N$, and
\eqref{eq:main_monotonicity} at every refinement level.
The recursion is uniformly bounded and monotone, and Taylor consistency
on smooth test functions gives the original operator in the limit.
To check the terminal condition, choose a smooth periodic $\varphi$ with
$\|\varphi-g\|_\infty\le\epsilon$. Bounded coefficients give a constant
$C_\varphi$, independent of $h,\tau,\nu$, such that
$\|\mathcal M_k[\varphi]-\varphi\|_\infty\le\tau C_\varphi$.
Monotonicity and constant-shift invariance then give
\[
 \varphi-\epsilon-C_\varphi(T-t_k)
 \le V_{h,\tau}^k\le
 \varphi+\epsilon+C_\varphi(T-t_k).
\]
Hence $\|V_{h,\tau}^k-g\|_\infty\le2\epsilon+C_\varphi(T-t_k)$
uniformly along the refinement. The Barles--Souganidis convergence result
\cite{barlessouganidis1991,guo2025policy} therefore gives
\[
 \max_k\|V_{h,\tau}^k-v(t_k)\|_\infty\to0.
\]
The first two terms
in~\eqref{eq:main_total_discrete} specify the required accuracy of the
neural approximation.

For the continuum estimate, retain continuous data and viscosity comparison.
Let $v$ be the continuous viscosity solution of the original equation,
and let $v_\nu$ be continuous viscosity solutions of the modified equation
with the same terminal data. Assume that $v_\nu$ are uniformly bounded
and that $\sup_{0<\nu\le1}\|v_\nu(t)-g\|_\infty\to0$ as $t\uparrow T$.
For the energy estimate, assume additionally the coefficient hypotheses
of Section~\ref{sec:stability} and that these viscosity solutions belong
to $W(0,T)$ and satisfy the variational equation.
Proposition~\ref{prop:main_wellposed} identifies them with the energy solutions above.
\begin{proposition}\label{prop:main_total}
Under these assumptions, $\|v_\nu-v\|_{L^\infty([0,T]\times\Omega)}\to0$.
For the iterates in Theorem~\ref{thm:main_history},
\begin{equation*}
\begin{aligned}
 \sup_t\|\widetilde v^n(t)-v(t)\|
 &\le\left[\frac73e^{\kappa_\rho T}
   \left(\sum_{m=1}^n\rho^{n-m}S_m+\rho^n\Phi_0\right)\right]^{1/2}\\
 &\quad+|\Omega|^{1/2}\|v_\nu-v\|_{L^\infty(Q_T)}.
\end{aligned}
\end{equation*}
\end{proposition}

\begin{proof}
The upper and lower half-relaxed limits of $v_\nu$ satisfy the original
viscosity inequalities: in a test function the added term
$\nu\Delta\varphi/2$ vanishes as $\nu\to0$. Uniform terminal continuity
gives the terminal data. Comparison identifies both limits with $v$,
and compactness of the space--time cylinder gives uniform convergence
\cite{barlesperthame1988,crandall1992guide}. The error estimate follows
by the triangle inequality and \eqref{eq:main_history}.
\end{proof}
Joint convergence follows when the weighted error sum, including its
$\lambda$-dependent prefactor, tends to zero. The discrete estimate uses
the stencil conditions~\eqref{eq:main_monotonicity}.

\subsection{Regularization bias for a smooth original solution}
For a smooth original solution $v$, we construct comparison barriers by
shifting $v$ by functions of time. This gives an explicit bound on the
regularization bias that is linear in $\nu$ and depends on
$\|\Delta v\|_\infty$.

\begin{corollary}\label{cor:main_smooth_bias}
Let $v\in C^{1,2}([0,T]\times\Omega)$ solve the original periodic equation
with terminal value $g$, and let $v_\nu$ be a continuous viscosity solution
of the equation with variances
$\Sigma_j+\nu$ and the same $f,c,\Sigma,g$. Assume viscosity comparison for
the modified equation. With $M=\|\Delta v\|_{L^\infty(Q_T)}$,
\begin{equation}\label{eq:main_smooth_bias}
 \|v_\nu(t)-v(t)\|_\infty\le\frac{\nu M}2(T-t),\qquad 0\le t\le T.
\end{equation}
\end{corollary}

\begin{proof}
Write $R_\nu[w]=w_t+H(t,x,\nabla w)+D_\nu w$. Since $v$ solves the
original equation, $R_\nu[v]=\nu\Delta v/2$.
For $w_\pm=v\pm\nu M(T-t)/2$,
\[
 R_\nu[w_+]=\tfrac\nu2(\Delta v-M)\le0,
 \qquad R_\nu[w_-]=\tfrac\nu2(\Delta v+M)\ge0.
\]
Both have terminal value $g$. Comparison for the backward equation gives
$w_-\le v_\nu\le w_+$, which proves the bound.
\end{proof}
\section{Numerical experiments}\label{sec:experiments}

We compare accuracy with direct PINN in five to fifty dimensions under
equal training time. Two-dimensional finite-difference method (FDM)
references separate neural error from viscosity bias, while comparisons
with direct FD training assess the cost of policy iteration.
Prescribed perturbations test the error-propagation estimate.
The refinement and dimension studies use three paired seeds and final
networks at the prescribed budget. Tables report means and sample
standard deviations (SDs).
Unless stated otherwise, neural experiments use three hidden layers of
width $64$ with sine activations, float64 arithmetic, batches of $256$,
and Adam with a cosine learning rate from $10^{-3}$ to $5\times10^{-5}$.
The schedule follows optimizer updates for fixed-update runs and elapsed
training time for fixed-time runs.
For evaluation points $Q=\{z_i\}_{i=1}^{N_Q}$, we use
\begin{equation*}
\begin{aligned}
 \operatorname{Err}_Q(w,v_{\rm ref})
 &=100\left(\frac{\sum_i|w(z_i)-v_{\rm ref}(z_i)|^2}
                    {\sum_i|v_{\rm ref}(z_i)|^2}\right)^{1/2},\\
 \operatorname{RMS}_Q(r)
 &=\left(\frac1{N_Q}\sum_i|r(z_i)|^2\right)^{1/2}.
\end{aligned}
\end{equation*}
The reference $v_{\rm ref}$ is either the original exact solution $v_\star$
or the FDM solution. FDM bias uses
$w=v_{\rm FDM}$ and $v_{\rm ref}=v_\star$. Gradient errors use the same formula
with vector gradients and Euclidean norms. Space--time errors pool the
evaluation points; fixed-time errors restrict the sample to one time.

\subsection{Degenerate HJB and HJI equations}\label{sec:fixed_family}
Let $x\in\mathbb T_{2\pi}^d$, $T=0.5$, $s=T-t$, and $S=\sum_jx_j$.
With cyclic indices, define
\begin{equation*}
\begin{aligned}
 g(x)&=\frac{0.3}{d}\sum_j\cos x_j+0.1\cos S,\\
 \psi(x)&=\frac1d\sum_j\sin x_j\sin x_{j+1}+0.4\sin S.
\end{aligned}
\end{equation*}
\begin{equation*}
 A(s)=s+0.05(1-\cos(2\pi s/T)),\qquad
 v_\star(t,x)=g(x)+A(T-t)\psi(x).
\end{equation*}
The physical variances alternate between $0.01$ and zero, so the
diffusion is rank deficient. The uncontrolled drift is
\[
 f_j^0(x)=0.15\sin x_{j+1}+0.05\sin(x_{j+2}-x_j).
\]
Dense orthonormal directions $q_k$ are generated by QR factorization
of a Gaussian matrix with seed $8421$; write $q_{jk}=(q_k)_j$.
For $d=2$, use the minimizing HJB Hamiltonian
\[
 H_0(x,p)=f^0(x)\cdot p-0.4\sum_{k=1}^{2}|q_k\cdot p|.
\]
For $d\in\{5,10,20,50\}$, add the maximizing term
$0.3\sum_{k=3}^{4}|q_k\cdot p|$ to obtain a nonconvex HJI Hamiltonian.
The full controlled drift is
$f(x,a,b)=f^0(x)+\sum_{k=1}^2a_kq_k+\sum_{k=3}^4b_kq_k$,
with $a_k\in[-0.4,0.4]$ and $b_k\in[-0.3,0.3]$; the last sum is absent
in the HJB case.
In every case, $H=H_0+c$ with the fixed source
\begin{equation*}
 c(t,x)=-\partial_t v_\star-H_0(x,\nabla v_\star)
                  -\frac12\sum_j\Sigma_j\partial_{x_jx_j}v_\star.
\end{equation*}
Thus $v_\star$ solves~\eqref{eq:main_original} with $f,c,\Sigma,g$
independent of $\nu$. The drift and control directions couple the
coordinates. The analytic derivatives and source identity were verified by AD.

\subsection{High-dimensional accuracy at equal training time}\label{sec:dimension_scaling}
For $d=5,10,20,50$, both methods use coordinate sine/cosine features
and the global features $\sin S$ and $\cos S$, which
appear in the manufactured solution $v_\star$.
There are $128d+8641$ trainable parameters, and each run has a
$180$-second training budget.

For FD-PINN-PI, set $\nu=0.02/d$ and
$h=0.75\nu/F_{\max}$, where
\[
 F_{\max}=\max_j\left\{0.2+0.4(|q_{j1}|+|q_{j2}|)
                              +0.3(|q_{j3}|+|q_{j4}|)\right\}.
\]
Choose a uniform time subdivision satisfying
$\tau\sum_j(\Sigma_j+\nu)/h^2\le0.8$.
For this family, $\|\Delta v_\star\|_\infty\le1.3+0.3d$; hence
\eqref{eq:main_smooth_bias} bounds the viscosity bias by
$\sup_t\|v_\nu-v_\star\|_\infty\le0.0015+0.0065/d$.
These choices satisfy~\eqref{eq:main_monotonicity}. The policy
interval is $200$ updates. Direct PINN uses $\nu=0$.
Both methods are evaluated against $v_\star$ on $16\,384$ independent
points in the full space--time domain.

\begin{table}[tbp]
\centering\small\begin{tabular}{rlrrr}
\toprule
$d$ & Method & Original error (\%) & Updates & Memory (MiB) \\
\midrule
5 & Direct PINN & $0.797 \pm 0.014$ & $18608 \pm 356$ & $75.9 \pm 0.0$ \\
5 & FD-PINN-PI & $0.538 \pm 0.055$ & $36117 \pm 663$ & $78.8 \pm 0.0$ \\
\midrule
10 & Direct PINN & $1.426 \pm 0.084$ & $13150 \pm 164$ & $82.1 \pm 0.0$ \\
10 & FD-PINN-PI & $1.008 \pm 0.068$ & $24000 \pm 566$ & $92.8 \pm 0.0$ \\
\midrule
20 & Direct PINN & $3.434 \pm 0.225$ & $7575 \pm 43$ & $99.4 \pm 0.0$ \\
20 & FD-PINN-PI & $1.814 \pm 0.074$ & $14750 \pm 25$ & $121.7 \pm 0.0$ \\
\midrule
50 & Direct PINN & $18.763 \pm 0.097$ & $3308 \pm 29$ & $165.4 \pm 0.0$ \\
50 & FD-PINN-PI & $15.230 \pm 0.754$ & $6283 \pm 14$ & $239.1 \pm 0.0$ \\
\bottomrule
\end{tabular}
\caption{Degenerate HJI errors after 180 seconds. Memory is peak PyTorch CUDA
allocation during training. Both methods use the features $\sin S,\cos S$.}
\label{tab:dimension_scaling}
\end{table}

FD-PINN-PI attains lower mean errors at every tested dimension and
completes roughly twice as many updates within the same training time,
with higher peak GPU allocation (Table~\ref{tab:dimension_scaling}).

\begin{figure}[tbp]
\centering\includegraphics[width=\textwidth]{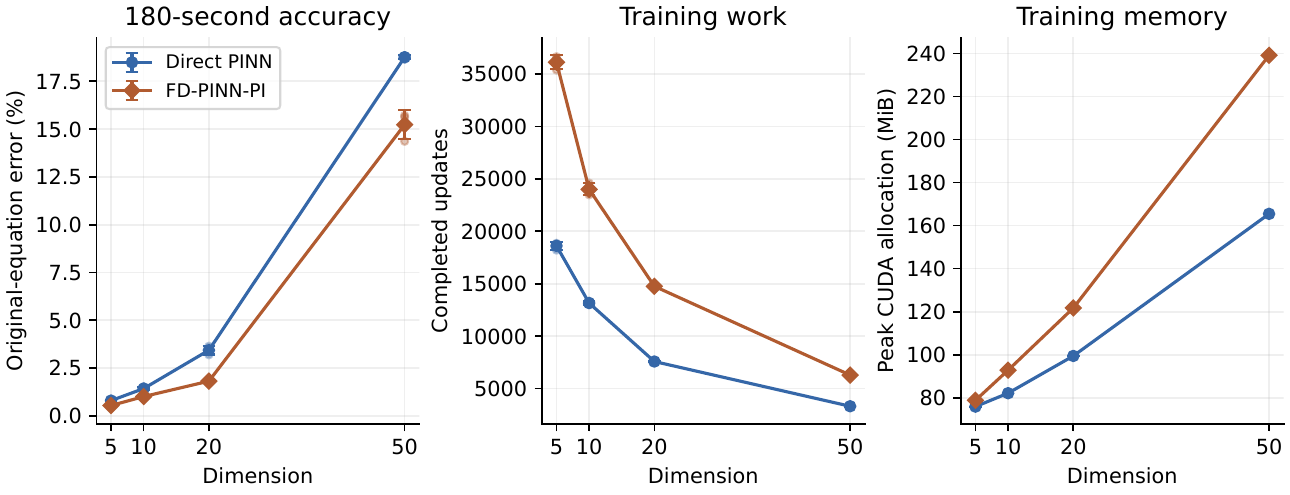}
\caption{Errors, optimizer updates, and peak GPU allocation after 180 seconds.}
\label{fig:dimension_scaling}
\end{figure}

\begin{figure}[tbp]
\centering\includegraphics[width=.94\textwidth]{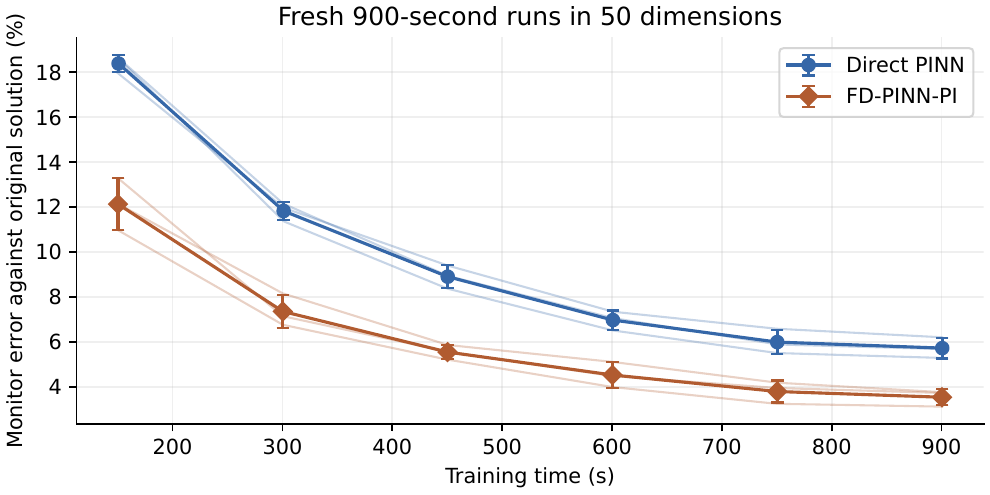}
\caption{Accuracy during the 900-second runs in fifty dimensions.
Thin lines show the three paired seeds; symbols and error bars give
means and one sample SD at the prescribed monitor times. Final table
entries use the independent test sample.}
\label{fig:main_extended_history}
\end{figure}

At $d=50$, separate 900-second runs use the same architecture, increments,
and paired initializations, with the cosine schedule scaled to this budget.
FD-PINN-PI reduces the mean relative error by about $40\%$ compared with
direct PINN (Table~\ref{tab:extended_fifty}), with lower final errors for
all three paired seeds (Figure~\ref{fig:main_extended_history}).

\begin{table}[tbp]
\centering\small\begin{tabular}{lrr}
\toprule
Method & Original error (\%) & Updates \\
\midrule
Direct PINN & $5.686 \pm 0.459$ & $16392 \pm 142$ \\
FD-PINN-PI & $3.386 \pm 0.409$ & $31408 \pm 14$ \\
\bottomrule
\end{tabular}
\caption{Fifty-dimensional errors after $900$ seconds, with three paired
initializations. The $180$-second results are in
Table~\ref{tab:dimension_scaling}; each budget uses separate training runs.}
\label{tab:extended_fifty}
\end{table}

\subsection{Low-dimensional accuracy and viscosity bias}\label{sec:fixed_target}
At $d=2$, the features are
$(2t/T-1,\sin x_1,\sin x_2,\cos x_1,\cos x_2)$.
Each network has $8769$ trainable parameters.
Both methods use batches of $512$ and $10\,000$ optimizer updates.
FD-PINN-PI refreshes the policy every $200$ updates and uses
\[
 (\nu,M)\in\{(0.08,64),(0.04,128),(0.02,256),(0.01,512)\},
 \qquad h=2\pi/M.
\]
The time increment is the largest step dividing $T$ and not exceeding
$h^2/4$. These increments satisfy~\eqref{eq:main_monotonicity}.
Direct PINN uses the original physical diffusion. Errors use $16\,384$
independent space--time points common to all runs.

FDM references use central differences, local-Lax--Friedrichs
stabilization, and the second-order strong-stability-preserving Runge--Kutta
method (SSPRK2) on $512^2$ and $1024^2$ periodic grids.
At these viscosities, the diffusion remains $\Sigma+\nu$.
Space--time interpolation evaluates both grids at the neural test points.
The maximum relative grid difference is \ReferenceGridDifference;
halving the fine-grid time step at $\nu=0.01$ gives \ReferenceTimeDifference.

\begin{table}[tbp]
\centering\footnotesize\setlength{\tabcolsep}{4pt}\begin{tabular}{lrrrrr}
\toprule
Method & $\nu$ & Original error (\%) & To FDM (\%) & FDM bias (\%) & Time (s) \\
\midrule
Direct PINN & $0$ & $0.252 \pm 0.016$ & --- & --- & $49.2 \pm 0.8$ \\
FD-PINN-PI & $0.08$ & $1.326 \pm 0.044$ & $0.282 \pm 0.017$ & $1.388$ & $42.3 \pm 0.3$ \\
FD-PINN-PI & $0.04$ & $0.707 \pm 0.047$ & $0.250 \pm 0.018$ & $0.697$ & $42.4 \pm 0.3$ \\
FD-PINN-PI & $0.02$ & $0.413 \pm 0.044$ & $0.252 \pm 0.018$ & $0.350$ & $42.3 \pm 0.1$ \\
FD-PINN-PI & $0.01$ & $0.297 \pm 0.033$ & $0.254 \pm 0.018$ & $0.175$ & $42.1 \pm 0.2$ \\
\bottomrule
\end{tabular}
\caption{Two-dimensional HJB errors (\%). Original error:
$\operatorname{Err}_Q(v_\theta,v_\star)$; To FDM:
$\operatorname{Err}_Q(v_\theta,v_{\rm FDM})$; FDM bias:
$\operatorname{Err}_Q(v_{\rm FDM},v_\star)$. The FDM solution uses the network's viscosity.}
\label{tab:fixed_target}
\end{table}

The FDM viscosity bias decreases at first order in $\nu$
(observed orders $0.99$--$1.00$), consistent with~\eqref{eq:main_smooth_bias}.
The neural error against the FDM reference remains near $0.25\%$,
so reducing viscosity improves accuracy against the original solution.
At the finest setting, direct PINN has a slightly lower error
(Table~\ref{tab:fixed_target}).

\begin{figure}[tbp]
\centering\includegraphics[width=.94\textwidth]{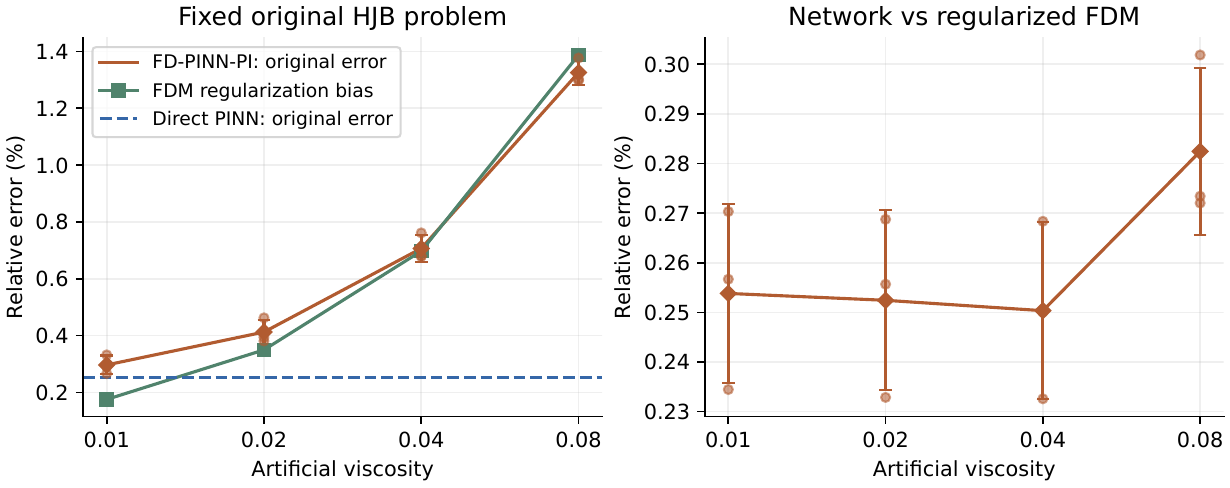}
\caption{Refinement with fixed original PDE data. Left: FD-PINN-PI error
against the original solution and the corresponding FDM regularization
bias. Right: network error against the positive-viscosity FDM reference.
Points show individual seeds; error bars indicate one sample SD.}
\label{fig:fixed_target}
\end{figure}

\subsection{Errors at fixed times and gradient accuracy}
\begin{figure}[tbp]
\centering\includegraphics[width=.96\textwidth]{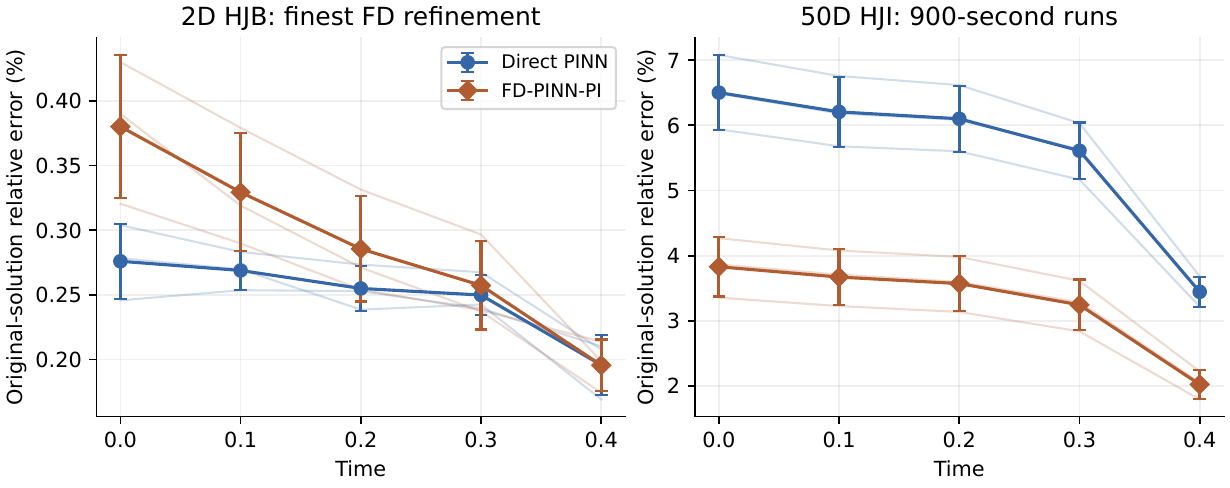}
\caption{Fixed-time errors against the original degenerate equations.
Thin curves show individual seeds; bars indicate one sample SD. The
fifty-dimensional FD-PINN-PI error is smaller at every evaluated time.}
\label{fig:main_time_slices}
\end{figure}

At $t=0,0.1,0.2,0.3,0.4$, we evaluate the finest two-dimensional runs
and the 900-second fifty-dimensional runs on $4096$ shared independent
spatial points (Figure~\ref{fig:main_time_slices}). The representation
enforces $v_\theta(T)=g$.
The relative gradient error is smaller for FD-PINN-PI in each tested
dimension (Table~\ref{tab:main_derivatives}). At $d=50$ and 180 seconds,
the errors are $4.876\%$ and $5.950\%$, respectively.
\begin{table}[tbp]
\centering\footnotesize\setlength{\tabcolsep}{3pt}
\begin{tabular}{rlrrrr}
\toprule
$d$ & Method & Gradient error (\%) & $R_0$ RMS & $R_\nu$ RMS & $R_h$ RMS \\
\midrule
5 & Direct PINN & $1.317 \pm 0.009$ & $0.0063 \pm 0.0001$ & $0.0063 \pm 0.0001$ & $0.0063 \pm 0.0001$ \\
5 & FD-PINN-PI & $0.838 \pm 0.061$ & $0.0045 \pm 0.0004$ & $0.0044 \pm 0.0004$ & $0.0044 \pm 0.0004$ \\
10 & Direct PINN & $1.707 \pm 0.086$ & $0.0094 \pm 0.0006$ & $0.0094 \pm 0.0006$ & $0.0094 \pm 0.0006$ \\
10 & FD-PINN-PI & $1.192 \pm 0.073$ & $0.0064 \pm 0.0005$ & $0.0063 \pm 0.0005$ & $0.0063 \pm 0.0005$ \\
20 & Direct PINN & $2.822 \pm 0.181$ & $0.0240 \pm 0.0020$ & $0.0240 \pm 0.0020$ & $0.0240 \pm 0.0020$ \\
20 & FD-PINN-PI & $1.626 \pm 0.037$ & $0.0115 \pm 0.0007$ & $0.0114 \pm 0.0007$ & $0.0114 \pm 0.0007$ \\
50 & Direct PINN & $5.950 \pm 0.160$ & $0.0801 \pm 0.0004$ & $0.0801 \pm 0.0004$ & $0.0801 \pm 0.0004$ \\
50 & FD-PINN-PI & $4.876 \pm 0.052$ & $0.0668 \pm 0.0030$ & $0.0668 \pm 0.0030$ & $0.0668 \pm 0.0030$ \\
\bottomrule
\end{tabular}
\caption{Gradient errors and residuals after 180 seconds.
Gradient errors are relative $L^2$ percentages. $R_0$ is the original
continuum residual, $R_\nu$ the method's regularized continuum residual,
and $R_h$ the nonlinear FD residual; residual columns are RMS values.
Direct PINN uses $\nu=0$. Entries are mean $\pm$ sample SD.}
\label{tab:main_derivatives}
\end{table}

\subsection{Air3D pursuit--evasion problem}\label{sec:air_main}
We use the Air3D relative-motion dynamics of
Mitchell et al.~\cite[Section III.B, Eqs.~(7)--(8)]{mitchell2005reachable},
setting both linear speeds to $0.5$ and both turn-rate bounds to $1.5$.
Let $(x,y)$ denote the relative planar position and $\psi$ the relative
heading angle. With momentum $p=(p_x,p_y,p_\psi)$, the Hamiltonian is
\begin{equation*}
 H=(-0.5+0.5\cos\psi)p_x+0.5\sin\psi\,p_y
       +1.5|yp_x-xp_y-p_\psi|-1.5|p_\psi|.
\end{equation*}
We solve the regularized terminal-value problem on $[-3,3]^2\times\mathbb S^1$
with zero physical diffusion, $T=0.5$, and artificial viscosity $\nu=0.2$.
All methods use $g_\epsilon=\sqrt{x^2+y^2+\epsilon^2}-0.5$ with
$\epsilon=0.1$ and the same planar Neumann penalty.
Five paired seeds use networks of width $128$ with features
\[
 (2t/T-1,x/3,y/3,\sin\psi,\cos\psi).
\]
The output is $g_\epsilon+(T-t)N_\theta$. Interior and boundary batches contain
$32\,768$ and $4096$ points. The PI policy is refreshed every $500$
updates, with FD increments $h=0.05$ and $\tau=0.0025$.
Training is in float32 with TF32 disabled. Errors use $16\,384$ points
with independently sampled coordinates in
$[-1.5,1.5]^2\times\mathbb S^1$ at $t=0,0.1,0.2,0.3,0.4$.
The relative error pools all five nonterminal slices.

\begin{table}[tbp]
\centering\small\begin{tabular}{lrrr}
\toprule
Method & Value error (\%) & Time (s) & Updates \\
\midrule
\multicolumn{4}{l}{\textit{30,000 updates}} \\
Direct PINN & $1.260 \pm 0.049$ & $807.7 \pm 0.7$ & $30\,000$ \\
FD-direct & $1.288 \pm 0.050$ & $390.7 \pm 0.5$ & $30\,000$ \\
AD-PINN-PI & $1.287 \pm 0.162$ & $832.5 \pm 5.7$ & $30\,000$ \\
FD-PINN-PI & $1.351 \pm 0.137$ & $476.1 \pm 4.2$ & $30\,000$ \\
\midrule
\multicolumn{4}{l}{\textit{600-second budget}} \\
Direct PINN & $1.334 \pm 0.159$ & $600.9 \pm 0.1$ & $22720 \pm 67$ \\
FD-direct & $1.174 \pm 0.129$ & $600.4 \pm 0.1$ & $46890 \pm 42$ \\
AD-PINN-PI & $1.442 \pm 0.157$ & $600.9 \pm 0.5$ & $21930 \pm 130$ \\
FD-PINN-PI & $1.261 \pm 0.118$ & $600.4 \pm 0.3$ & $38200 \pm 117$ \\
\bottomrule
\end{tabular}
\caption{Four-method Air3D value errors against the independently refined FDM
reference at common $\nu=0.2$. Entries are mean $\pm$ sample SD over
five paired seeds, for fixed updates and separate fixed-time runs.}
\label{tab:main_air}
\end{table}
At fixed updates, FD-PINN-PI takes about $41\%$ less training time than
direct PINN, with a slightly higher mean error. At equal training time,
its mean error is lower (Table~\ref{tab:main_air}). Among all four
methods, FD-direct attains the lowest mean error at this budget, $1.174\%$.

Reference differences are $0.00446\%$ under grid refinement and
$0.00268\%$ under time-step halving at $t=0$; domain enlargement changes
the pooled evaluation values by $0.00187\%$.
At the same terminal smoothing, the difference from the reference with
zero artificial viscosity is $5.669\%$.
\ref{subsec:air3d_revision} gives the refinement and boundary diagnostics.

\subsection{Policy-update interval}\label{sec:main_policy}
We compare FD-PINN-PI with FD-direct for the twenty-dimensional HJI
problem, whose controls have closed-form updates. The PDE, viscosity,
architecture, and 180-second budget are fixed. FD-direct attains a lower
mean error than each tested PI interval (Table~\ref{tab:main_policy}).

Table~\ref{tab:main_policy} reports the nonlinear residual and the gradient lag
$\operatorname{RMS}_Q(\nabla_h v_\theta-p)$ on $2048$ fresh points,
where $p$ is the frozen network's FD gradient. This gradient difference
is the second term in~\eqref{eq:main_discrete_mismatch}.
\begin{table}[tbp]
\centering\footnotesize\setlength{\tabcolsep}{3pt}
\begin{tabular}{lrrrr}
\toprule
Method & Original error (\%) & Updates & Nonlinear residual & Gradient lag \\
\midrule
FD-direct & $1.583 \pm 0.064$ & $17217 \pm 38$ & $0.00970 \pm 0.00053$ & $0.00000 \pm 0.00000$ \\
$K=50$ & $1.781 \pm 0.058$ & $14733 \pm 38$ & $0.01110 \pm 0.00012$ & $0.00130 \pm 0.00009$ \\
$K=200$ & $1.814 \pm 0.074$ & $14750 \pm 25$ & $0.01134 \pm 0.00048$ & $0.00155 \pm 0.00015$ \\
$K=1000$ & $1.754 \pm 0.068$ & $14750 \pm 43$ & $0.01099 \pm 0.00052$ & $0.00145 \pm 0.00019$ \\
\bottomrule
\end{tabular}
\caption{Twenty-dimensional comparison with the same original PDE,
positive-viscosity stencil, and 180-second budget. Entries are mean
$\pm$ sample SD over three paired seeds. Gradient lag is zero by
definition for direct nonlinear FD training.}
\label{tab:main_policy}
\end{table}

\subsection{Isaacs equation with finite control sets}\label{sec:finite_controls}
To examine the cost of control optimization, we consider a game in which
each policy update searches a table of control pairs.
Retain the domain, rank-deficient diffusion, and exact solution
$v_\star$ of Section~\ref{sec:fixed_family} at $d=20$, and replace $H_0$ by
\[
 H_M(x,p)=\min_{0\leq i<M}\max_{0\leq j<M}
 \{f_{ij}(x)\cdot p+\ell_{ij}(x)\},\qquad M=512.
\]
Replacing $v$ by $-v$ and interchanging the players gives the max--min
convention of~\eqref{eq:main_nested}.
Let $U\in\mathbb R^{20\times6}$ have orthonormal columns and set
$\alpha_i=2\pi i/M$, $\beta_j=2\pi j/M$. Define
\[
\begin{aligned}
 z_{ij}&=(-0.4\cos\alpha_i,-0.4\sin\alpha_i,
          0.3\cos\beta_j,0.3\sin\beta_j,r_{ij,1},r_{ij,2}),\\
 f_{ij}(x)&=f^0(x)+Uz_{ij},\\
 \ell_{ij}(x)&=0.02\sin x_1\,\eta_{ij,1}
                      +0.015\cos x_2\,\eta_{ij,2}.
\end{aligned}
\]
For $k=1,2$, draw $r_{ij,k}$ uniformly from $[-0.15,0.15]$ and
$\eta_{ij,k}$ uniformly from $[-1,1]$, independently.
The matrix $U$ is obtained by Gaussian QR factorization; the matrix and
tables are fixed with seed $8421$. The coupled coefficients require
a nested search over the control pairs. The source
$c=-\partial_t v_\star-H_M(x,\nabla v_\star)
-\frac12\sum_k\Sigma_k\partial_{x_kx_k}v_\star$
keeps the original degenerate solution fixed.

All methods use $\nu=0.005$, $h=0.00601508$, and
$\tau=0.000144725$, satisfying~\eqref{eq:main_monotonicity}.
We retain the periodic and terminal representation of the dimension
study and use $10\,000$ optimizer updates.
FD-direct and FD-PINN-PI share initializations and reuse identical
batches for $K=5$ updates. FD-direct computes the exact nested optimum
at every update; PI computes it once and retains the selected controls.
Both differentiate the selected affine branch. We also run FD-direct
with fresh samples at every update.

After pilot runs with seed zero, the configuration is fixed for seeds
$1,2,3$. On a shared RTX~5070, we alternate synchronized training blocks
and vary the method order. Times include sampling and source evaluation,
performed once per retained batch or at every update with fresh sampling.
We check $2048$ monitor points every $500$ updates and validate the first
checkpoint below $5\%$ error on $16\,384$ independent test points.

\begin{table}[!ht]
\centering\footnotesize\setlength{\tabcolsep}{4pt}
\begin{tabular}{lrrr}
\toprule
Method & Time to $5\%$ (s) & Final error (\%) & Total time (s) \\
\midrule
FD-direct, reused samples & $143.7 \pm 6.6$ & $3.253 \pm 0.126$ & $226.5 \pm 1.3$ \\
FD-PINN-PI, cached policy & $93.9 \pm 4.2$ & $3.251 \pm 0.116$ & $148.0 \pm 0.9$ \\
FD-direct, fresh samples & $124.2 \pm 1.8$ & $2.690 \pm 0.143$ & $310.4 \pm 1.3$ \\
\bottomrule
\end{tabular}
\caption{Twenty-dimensional game with $512$ controls per player: mean $\pm$ sample SD
over three paired seeds. Final errors and total times use $10\,000$
updates. All nine threshold checkpoints satisfy $5\%$ error on the
independent test set.}\label{tab:finite_controls}
\end{table}
\begin{figure}[!ht]
\centering\includegraphics[width=.98\textwidth]{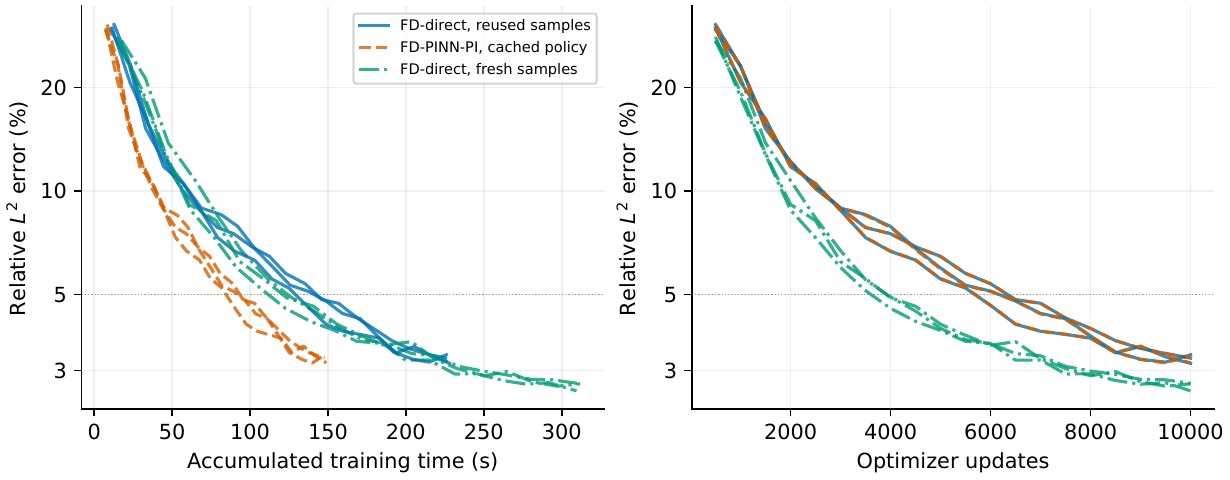}
\caption{Game with finite control sets: monitor error against accumulated training
time (left) and optimizer updates (right). Each method has three seed
curves; the horizontal line marks $5\%$ error.}
\label{fig:finite_controls}
\end{figure}

With identical samples reused for $K=5$ updates, PI reduces the mean time to $5\%$ error from
$143.7$ to $93.9$ seconds, or $34.6\%$, while the errors after
$10\,000$ updates are nearly identical
(Table~\ref{tab:finite_controls}, Figure~\ref{fig:finite_controls}).
The number of states requiring a policy search falls from $2\,560\,000$
to $512\,000$. After $10\,000$ updates, FD-direct with fresh samples has
a lower final error, $2.690\%$, while PI requires $24.3\%$ less time to
reach $5\%$ error. At approximately $140$ seconds, the mean relative error
on the common monitoring set is $3.27\%$ for PI and $4.42\%$ for
FD-direct with fresh samples.
The comparison with identical samples attributes the saving to fewer
control optimizations during policy evaluation.

\subsection{Prescribed policy-iteration errors}\label{sec:main_defects}
To test how policy-evaluation errors accumulate in~\eqref{eq:main_history},
we prescribe perturbations on $\mathbb T_{2\pi}^2$ with $T=0.5$,
$\Sigma=(0.01,0)$, $\nu=0.25$, and Hamiltonian
$H=-0.12|p_1|+0.08|p_2|+c$. The exact solution is
$v_\star=g+(T-t)\psi$, where
\[
\begin{aligned}
 g&=\cos x_1+\tfrac12\sin x_2+\tfrac14\cos(x_1+x_2),\\
 \psi&=0.4\sin x_1\sin x_2+0.2\cos(2x_1).
\end{aligned}
\]
The source is manufactured for this fixed modified equation. Nested
improvement gives the discontinuous controls
$\alpha=-0.12\operatorname{sign}(p_1)$ and
$\beta=0.08\operatorname{sign}(p_2)$. Policy evaluation uses Fourier
spatial discretization and RK4 time stepping.

Momentum, control optimization, and terminal perturbations are activated separately
and together. Their amplitudes are persistent, decrease as $0.55^{n-1}$,
or occur only at iteration two. Each computed iterate is reconstructed
as a trigonometric polynomial in space and a piecewise-linear function
in time; all four continuum defect norms are recomputed on independent
shifted quadratures. Let $\overline E_n=E_n/|\Omega|$ be the energy error
per unit volume. With $S_m$ and $\Phi_0$ as in~\eqref{eq:main_history}, define
\begin{equation}\label{eq:main_diagnostic_history}
\begin{aligned}
 B_n&=\frac{7e^{\kappa_\rho T}}{3|\Omega|}
       \left(\sum_{m=1}^n\rho^{n-m}S_m+\rho^n\Phi_0\right),\\
 B_n^r&=\frac{7e^{\kappa_\rho T}}{3|\Omega|}
       \left(\sum_{m=1}^n\rho^{n-m}\|r_m\|_{L^2(Q_T)}^2
                                     +\rho^n\Phi_0\right).
\end{aligned}
\end{equation}
Here $\rho=0.5$, $\lambda=0.25$, $F^2=0.12^2+0.08^2$, and
$K_\Sigma=0$.

\begin{figure}[tbp]
\centering\includegraphics[width=.98\textwidth]{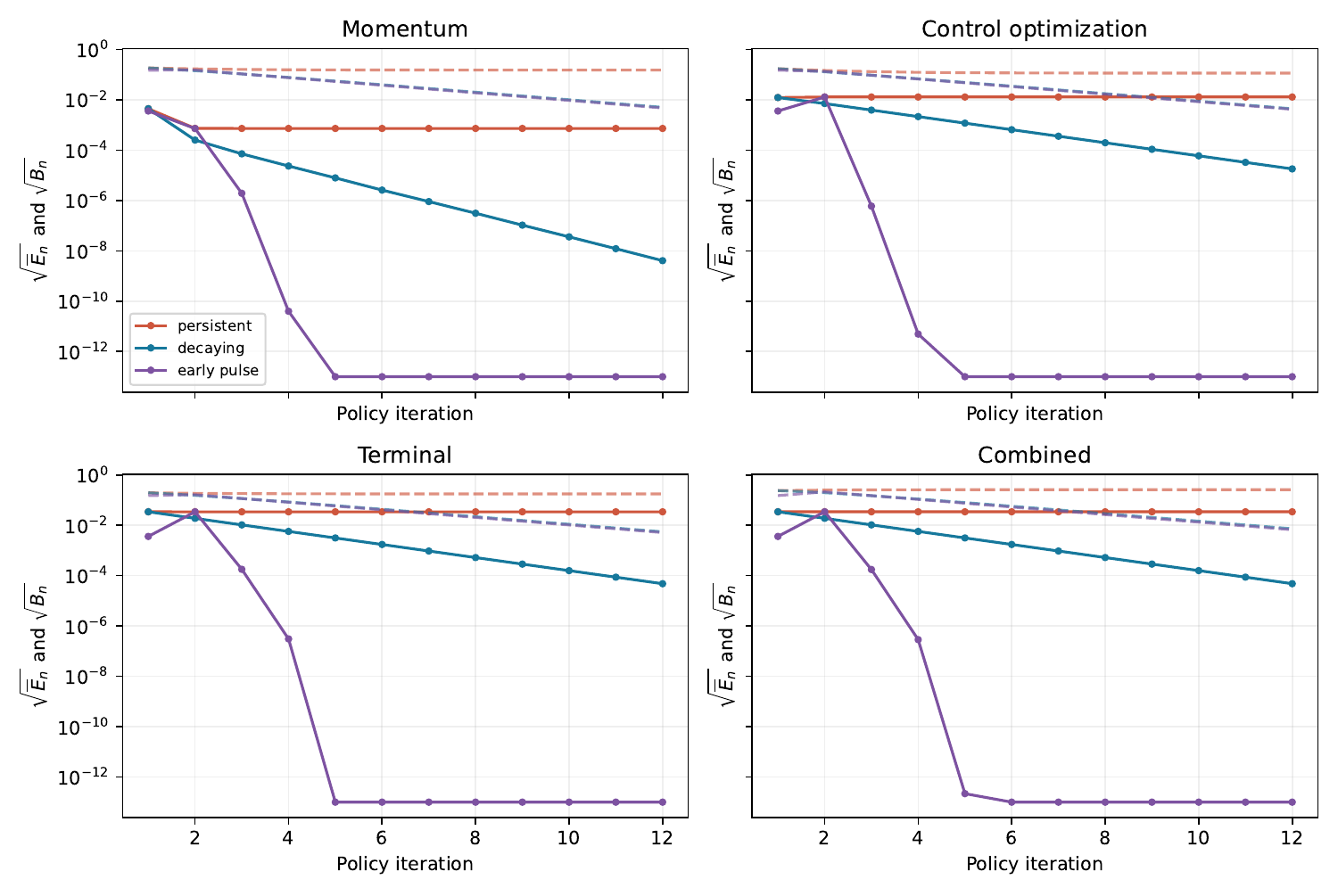}
\caption{Prescribed PI errors with discontinuous feedbacks.
Solid curves: $\sqrt{\overline E_n}$; dashed curves: $\sqrt{B_n}$.
Perturbation amplitudes are constant, proportional to $0.55^{n-1}$,
or nonzero only at $n=2$.}
\label{fig:main_defects}
\end{figure}
\begin{table}[tbp]
\centering\small
\begin{tabular}{lrrr}
\toprule
Perturbed term & $\sqrt{\overline E_{12}}$ & $\overline E_{12}/B_{12}$ & $\overline E_{12}/B^r_{12}$\\
\midrule
Momentum & $0.00073$ & $0.00002$ & $0.011$\\
Control optimization & $0.01293$ & $0.01268$ & $15.054$\\
Terminal & $0.03419$ & $0.03795$ & $98.111$\\
Combined & $0.03425$ & $0.01738$ & $35.063$\\
\bottomrule
\end{tabular}
\caption{Constant perturbations after twelve PI steps.
$\overline E_{12}/B_{12}<1$ in all cases; $\overline E_{12}/B^r_{12}>1$ for the control optimization,
terminal, and combined perturbations.}
\label{tab:main_defects}
\end{table}
Across all twelve runs, $\overline E_n/B_n<0.044$. The perturbations and
spatial, temporal, and quadrature refinements are specified in
\ref{subsec:periodic_defects}.

\section{Conclusion}\label{sec:conclusion}

FD-PINN-PI extends neural policy evaluation to high-dimensional degenerate
Hamilton--Jacobi equations through viscosity-based semi-discrete policy
iteration. Local finite differences define the control update and evaluation
residual, while the network represents the value function without a global
tensor grid. The analysis bounds the propagation of errors for measurable
feedbacks and relates the discrete approximation to the original equation.

On the coupled HJI benchmarks, the method gives lower mean errors than
direct PINN from five to fifty dimensions at equal training time, with
about $40\%$ lower error at fifty dimensions and $900$ seconds.
Two-dimensional references verify the decrease in viscosity bias.
In the tested problems with closed-form controls, FD-direct attains
lower mean errors.
For the twenty-dimensional game with finite control sets, reusing policies
reduces the measured time to $5\%$ error by $24.3\%$ relative to FD-direct
with fresh samples at every update.
\appendix
\renewcommand{\thetheorem}{\Alph{section}.\arabic{theorem}}
\section{Well-posedness of neural policy-evaluation targets}\label{app:wellposed}
Write
$f_{\alpha,\beta}=f(t,x,\alpha,\beta)$ and
$c_{\alpha,\beta}=c(t,x,\alpha,\beta)$.
\begin{proof}[Proof of Proposition~\ref{prop:main_wellposed}]
Reverse time by $s=T-t$, so that~\eqref{eq:main_wellposed_frozen} becomes an initial-value problem. Using~\eqref{eq:main_divergence}, the spatial operator is realized by the time-dependent bilinear form on $H^1(\Omega)$,
\[
\begin{aligned}
\mathfrak a_s(u,\varphi)
&:=\frac12\sum_{j=1}^d\int_\Omega(\Sigma_j+\nu)\partial_{x_j}u\,\partial_{x_j}\varphi\,dx\\
&\quad+\frac12\sum_{j=1}^d\int_\Omega(\partial_{x_j}\Sigma_j)(\partial_{x_j}u)\varphi\,dx
-\int_\Omega(f_{\alpha,\beta}\cdot\nabla u)\varphi\,dx .
\end{aligned}
\]
Its coefficients are measurable in $s$ and bounded by the stated coefficient assumptions. The form is bounded, $|\mathfrak a_s(u,\varphi)|\le\bigl(\tfrac{\Lambda}{2}+\tfrac{K_\Sigma}{2}+F\bigr)\|u\|_{H^1}\|\varphi\|_{H^1}$, and G\aa rding-coercive:
\[
 \mathfrak a_s(u,u)
 \ge
 \frac{\lambda}{2}\|\nabla u\|^2
 -\Bigl(F+\frac{K_\Sigma}{2}\Bigr)\|\nabla u\|\|u\|
 \ge
 \frac{\lambda}{4}\|\nabla u\|^2
 -\frac{(F+K_\Sigma/2)^2}{\lambda}\|u\|^2 .
\]
Adding $\gamma\|u\|^2$, with
$\gamma=1+(F+K_\Sigma/2)^2/\lambda$, gives coercivity on $H^1(\Omega)$.
Lions' theorem for non-autonomous forms therefore yields a unique solution
in $W(0,T)$ for every $L^2(Q_T)$ source and $L^2(\Omega)$ initial value
\cite{lions1961,dautraylions1992}. Reversing time gives the solution of~\eqref{eq:main_wellposed_frozen}.

Let $X:=L^2(0,T;H^1(\Omega))$. The Hamiltonian is measurable in $(t,x)$ as a $\sup$--$\inf$ of
Carath\'eodory functions over compact separable control sets, and it is
$F$-Lipschitz in $p$ by \eqref{eq:L_H_lipschitz}. Consequently,
$(t,x)\mapsto H(t,x,\nabla w(t,x))$ is measurable for every $w\in X$.
It belongs to $L^2(Q_T)$ by $|H(t,x,p)|\le C_c+F|p|$, where
$C_c:=\esssup_{(t,x)\in Q_T}\sup_{a,b}|c(t,x,a,b)|$.
The linear existence result, with zero drift and the $L^2$ source
$-H(\cdot,\cdot,\nabla w)$, defines $\mathcal T[w]\in W(0,T)$ as the
unique variational solution of
\[
 \partial_tv+\frac12\sum_{j=1}^d(\Sigma_j+\nu)\partial_{x_jx_j}^2v
 =-H(t,x,\nabla w),
 \qquad v(T)=g .
\]
Fixed points of $\mathcal T$ are exactly the variational solutions of~\eqref{eq:main_target}. For $w_1,w_2\in X$, the difference $z:=\mathcal T[w_1]-\mathcal T[w_2]$ solves the same linear problem with terminal value $0$ and source $G:=-\bigl(H(\cdot,\cdot,\nabla w_1)-H(\cdot,\cdot,\nabla w_2)\bigr)$, so that $\|G(t)\|\le F\|\nabla(w_1-w_2)(t)\|$. Testing with $z$ as in \eqref{eq:frozen_energy_estimate} gives, for every $\mu>0$,
\[
 -\frac{d}{dt}\|z\|^2
 +\frac{3\lambda}{4}\|\nabla z\|^2
 \le
 \Bigl(\frac{K_\Sigma^2}{\lambda}+\mu\Bigr)\|z\|^2
 +\frac1\mu\|G\|^2 .
\]
Multiplying by $\omega_\mu(t):=e^{-(K_\Sigma^2/\lambda+\mu)(T-t)}$ and integrating over $(t,T)$ with $z(T)=0$,
\[
 \omega_\mu(t)\|z(t)\|^2
 +\frac{3\lambda}{4}\int_t^T\omega_\mu\|\nabla z\|^2\,ds
 \le
 \frac{F^2}{\mu}\int_t^T\omega_\mu\|\nabla(w_1-w_2)\|^2\,ds .
\]
Taking the supremum after discarding the gradient term, and taking $t=0$ after discarding the value term, bounds the two quantities separately. Consequently, in the norm $|||w|||^2:=\int_0^T\omega_\mu\bigl(\|w\|^2+\|\nabla w\|^2\bigr)dt$, which is equivalent to the norm of $X$ because $e^{-(K_\Sigma^2/\lambda+\mu)T}\le\omega_\mu\le1$,
\[
 |||\mathcal T[w_1]-\mathcal T[w_2]|||^2
 \le
 \frac{F^2}{\mu}\Bigl(\frac{4}{3\lambda}+T\Bigr)\,|||w_1-w_2|||^2 .
\]
Choosing $\mu>F^2\bigl(\tfrac{4}{3\lambda}+T\bigr)$ makes $\mathcal T$ a strict contraction on the complete space $(X,|||\cdot|||)$. The unique fixed point $v_{\nu}=\mathcal T[v_{\nu}]$ belongs to
$W(0,T)$ and is the unique variational solution of~\eqref{eq:main_target}.
\end{proof}
\section{Independent residual validation}\label{app:validation}
\begin{lemma}\label{lem:analysis_sampling}
Fix an iteration $n$. Let $\mathcal G_n$ be the $\sigma$-algebra generated by all randomness determining the trained surrogate, the frozen feedback pair, and any model selection used to obtain them. Let $Z_1,\ldots,Z_N$ be i.i.d.\ uniform points on $Q_T$, independent of $\mathcal G_n$. Suppose a finite $\mathcal G_n$-measurable envelope $\overline B_n$ is known such that
\[
 |q_n(t,x)|\le\overline B_n
 \quad\text{for a.e. }(t,x)\in Q_T.
\]
Set
\[
 \widehat{\mathcal J}_{n,N}^{\mathrm{val}}
 :=\frac1N\sum_{i=1}^N|q_n(Z_i)|^2.
\]
Then, for every $\delta\in(0,1)$, with probability at least $1-\delta$,
\begin{equation}\label{eq:validation_bound}
 \|q_n\|_{L^2(Q_T)}^2
 \le
 |Q_T|\left(
 \widehat{\mathcal J}_{n,N}^{\mathrm{val}}
 +\overline B_n^2\sqrt{\frac{\log(1/\delta)}{2N}}
 \right).
\end{equation}
If $\widetilde v^n\in C^2([-\tau,T];C^4(\Omega))$, then also
\begin{equation}\label{eq:validated_continuum_residual}
 \|r_n\|_{L^2(Q_T)}^2
 \le
 2|Q_T|\left(
 \widehat{\mathcal J}_{n,N}^{\mathrm{val}}
 +\overline B_n^2\sqrt{\frac{\log(1/\delta)}{2N}}
 +\mathfrak c_{\tau,h}[\widetilde v^n]^2
 \right).
\end{equation}
\end{lemma}

\begin{proof}
Conditionally on $\mathcal G_n$, the variables $|q_n(Z_i)|^2$ are i.i.d. and lie in $[0,\overline B_n^2]$. Their common expectation is $|Q_T|^{-1}\|q_n\|_{L^2(Q_T)}^2$. Hoeffding's inequality yields~\eqref{eq:validation_bound}. Combining this with~\eqref{eq:residual_consistency_l2} and $(a+b)^2\le2a^2+2b^2$ gives~\eqref{eq:validated_continuum_residual}; averaging the conditional probability gives the stated unconditional bound.
\end{proof}

The bound requires a global residual envelope and validation points drawn
independently after network selection.
For finitely many policy stages, applying~\eqref{eq:validation_bound} with failure
probabilities $\delta_m$ satisfying $\sum_m\delta_m\le\delta$ gives
simultaneous bounds by the union bound, provided each validation sample
is independent of the corresponding training and selection history.
Temporal consistency cannot be omitted: the spatially constant function
$w(t,x)=\sin(2\pi(t-T)/\tau)$ has zero time-quotient residual and zero terminal
value for the zero-data problem, while its continuum time derivative is
nonzero. At $t_k=k\tau$ it vanishes, consistently with~\eqref{eq:main_discrete_certificate}.
\counterwithin{table}{section}
\counterwithin{figure}{section}
\renewcommand{\thetable}{\Alph{section}.\arabic{table}}
\renewcommand{\thefigure}{\Alph{section}.\arabic{figure}}
\section{HJB and HJI numerical data}\label{sec:priority_experiments}
\subsection{Sampling, implementation, and coefficient checks}\label{subsec:priority_protocol}
Training, monitoring, and final evaluation use independent random
generators, with common evaluation points across methods.

The implementation verifies analytic derivatives against AD,
$R_\nu[v_\star]=\nu\Delta v_\star/2$ and second-order spatial/first-order
temporal consistency along $(h,\tau)\mapsto(h/2,\tau/4)$.
Table~\ref{tab:priority_coefficients} gives the FD parameters;
the spatial margin is $\min_j(\Sigma_j+\nu-hF_j)$. The coordinate drift
bound is $F_j=0.2+0.4\sum_{k=1}^2|q_{jk}|+0.3\sum_{k=3}^4|q_{jk}|$,
with the final sum omitted in the two-dimensional HJB case.
\begin{table}[!htbp]
\centering\small\begin{tabular}{rrrrrr}
\toprule
$d$ & $\nu$ & $h$ & $\tau$ & Spatial margin & Time CFL \\
\midrule
2 & 0.08 & 0.0981748 & 0.00240385 & 0.0207881 & 0.042 \\
2 & 0.04 & 0.0490874 & 0.000601685 & 0.0103941 & 0.022 \\
2 & 0.02 & 0.0245437 & 0.000150557 & 0.00519704 & 0.012 \\
2 & 0.01 & 0.0122718 & 3.76478e-05 & 0.00259852 & 0.007 \\
5 & 0.004 & 0.00382415 & 0.000233973 & 0.00115407 & 0.800 \\
10 & 0.002 & 0.00213698 & 5.21866e-05 & 0.0005 & 0.800 \\
20 & 0.001 & 0.00114658 & 8.76424e-06 & 0.00025 & 0.800 \\
50 & 0.0004 & 0.00055304 & 9.06229e-07 & 0.0001 & 0.800 \\
\bottomrule
\end{tabular}
\caption{FD increments and coefficient checks. Time CFL is
$\tau\sum_j(\Sigma_j+\nu)/h^2$. All configurations satisfy
\eqref{eq:main_monotonicity}.}\label{tab:priority_coefficients}
\end{table}

\subsection{Independent FDM references}
Let $h_{\rm ref}=2\pi/M_{\rm ref}$ denote the reference-grid spacing.
The periodic reference uses central differences with variance
$a_j=\max\{\Sigma_j+\nu,h_{\rm ref}F_j\}$ and SSPRK2 time stepping. This stencil
uses the monotone comparison structure of~\cite{barlessouganidis1991}.
The independent analytic source is checked against the neural source.
Grids $512^2$ and $1024^2$ use
$\Delta t\le\min\{0.001,0.8h_{\rm ref}^2/\sum_j a_j\}$,
adjusted to divide $T$.
A second fine solve at $\nu=0.01$ halves the time-step cap.
Positive-viscosity references require no added stabilization; at zero
artificial viscosity the stated local-Lax--Friedrichs variance remains
in degenerate directions. Periodic bilinear interpolation and linear
time interpolation evaluate the references at common neural test points.
\begin{table}[!htbp]
\centering\small\begin{tabular}{rrrrr}
\toprule
$\nu$ & Grid & $\Delta t$ & Steps & Extra variance \\
\midrule
0.08 & 512 & 0.000708215 & 706 & 0 \\
0.08 & 1024 & 0.000177117 & 2823 & 0 \\
0.04 & 512 & 0.001 & 500 & 0 \\
0.04 & 1024 & 0.000334448 & 1495 & 0 \\
0.02 & 512 & 0.001 & 500 & 0 \\
0.02 & 1024 & 0.000601685 & 831 & 0 \\
0.01 & 512 & 0.001 & 500 & 0 \\
0.01 & 1024 & 0.0005 & 1000 & 0 \\
0.01 & 1024 & 0.001 & 500 & 0 \\
0.00 & 512 & 0.001 & 500 & 0.00740148 \\
0.00 & 1024 & 0.001 & 500 & 0.00370074 \\
\bottomrule
\end{tabular}
\caption{Independent reference discretizations. Extra variance is
$\max_j(a_j-\Sigma_j-\nu)$.}\label{tab:priority_reference_details}
\end{table}

\subsection{Fixed-time evaluations}\label{subsec:time_slices}
All $45$ primary and longer-budget checkpoints are evaluated at
$t=0,0.1,0.2,0.3,0.4$ using $4096$ shared spatial points drawn
independently of training, monitor, and space--time tests.
The two-dimensional points also query the FDM references.
Figure~\ref{fig:main_time_slices} plots the profiles;
Tables~\ref{tab:initial_fixed_target} and
\ref{tab:initial_dimension_study} give all initial-time results.
\begin{table}[!htbp]
\centering\small\begin{tabular}{lrrrr}
\toprule
Method & $\nu$ & Original error (\%) & To FDM (\%) & FDM bias (\%) \\
\midrule
Direct PINN & 0.00 & $0.276 \pm 0.029$ & --- & --- \\
FD-PINN-PI & 0.08 & $2.108 \pm 0.047$ & $0.260 \pm 0.057$ & $2.125$ \\
FD-PINN-PI & 0.04 & $1.089 \pm 0.056$ & $0.268 \pm 0.054$ & $1.070$ \\
FD-PINN-PI & 0.02 & $0.595 \pm 0.059$ & $0.274 \pm 0.053$ & $0.537$ \\
FD-PINN-PI & 0.01 & $0.380 \pm 0.055$ & $0.277 \pm 0.052$ & $0.269$ \\
\bottomrule
\end{tabular}
\caption{Initial-time two-dimensional errors. Original error is against
$v_\star$; the remaining columns separate neural--FDM difference and
regularization bias. Neural entries are mean $\pm$ sample SD.}
\label{tab:initial_fixed_target}
\end{table}
\begin{table}[!htbp]
\centering\small\begin{tabular}{rrlrr}
\toprule
$d$ & Budget (s) & Method & Original error (\%) & Absolute RMS \\
\midrule
5 & 180 & Direct PINN & $0.855 \pm 0.001$ & $0.00185 \pm 0.00000$ \\
5 & 180 & FD-PINN-PI & $0.602 \pm 0.022$ & $0.00130 \pm 0.00005$ \\
\midrule
10 & 180 & Direct PINN & $1.616 \pm 0.112$ & $0.00306 \pm 0.00021$ \\
10 & 180 & FD-PINN-PI & $1.207 \pm 0.080$ & $0.00228 \pm 0.00015$ \\
\midrule
20 & 180 & Direct PINN & $3.901 \pm 0.222$ & $0.00682 \pm 0.00039$ \\
20 & 180 & FD-PINN-PI & $2.072 \pm 0.088$ & $0.00362 \pm 0.00015$ \\
\midrule
50 & 180 & Direct PINN & $21.471 \pm 0.114$ & $0.03502 \pm 0.00019$ \\
50 & 180 & FD-PINN-PI & $17.109 \pm 0.739$ & $0.02791 \pm 0.00120$ \\
\midrule
50 & 900 & Direct PINN & $6.502 \pm 0.570$ & $0.01061 \pm 0.00093$ \\
50 & 900 & FD-PINN-PI & $3.830 \pm 0.457$ & $0.00625 \pm 0.00075$ \\
\bottomrule
\end{tabular}
\caption{Initial-time accuracy in the dimension and budget comparisons.
Relative errors are percentages; absolute RMS uses the same spatial points.}
\label{tab:initial_dimension_study}
\end{table}

\subsection{Policy-update comparison}\label{subsec:policy_intervals}
\begin{figure}[!htbp]
\centering\includegraphics[width=.94\textwidth]{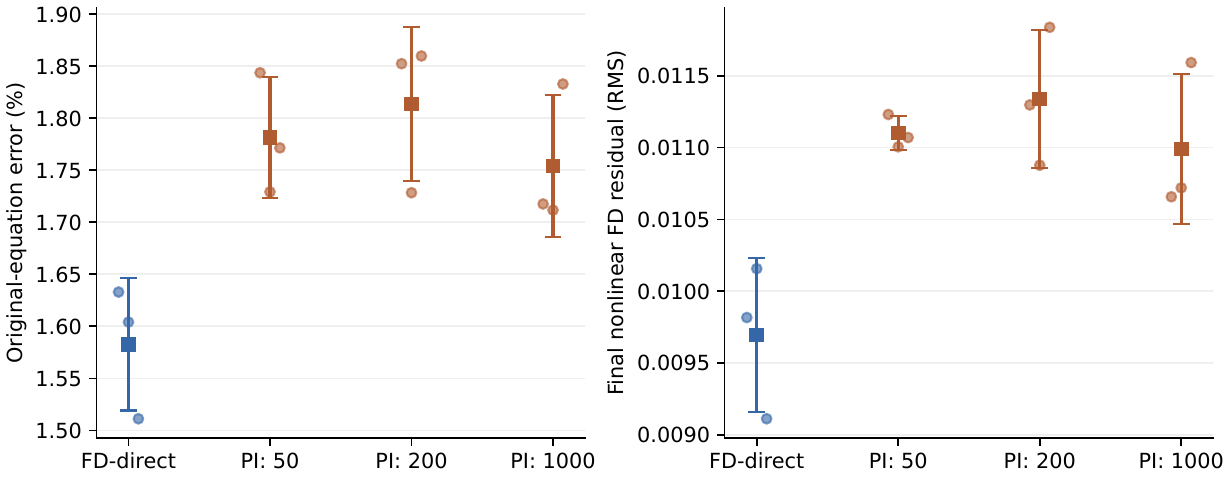}
\caption{Individual twenty-dimensional errors and nonlinear FD residuals
at the final 180-second budget. Symbols and bars show mean and sample SD.}
\label{fig:policy_interval}
\end{figure}

The twenty-dimensional comparison retains the three primary $K=200$
runs and adds three paired seeds for FD-direct, $K=50$, and $K=1000$.
Final policies use $2048$ independent validation points to check
the control optimization error and~\eqref{eq:main_discrete_mismatch}.
When the frozen and current networks coincide, direct and PI residuals
and their parameter gradients agree. Figure~\ref{fig:policy_interval}
shows the individual runs summarized in Table~\ref{tab:main_policy}.

\section{Air3D numerical data}\label{subsec:air3d_revision}
The four methods in Table~\ref{tab:main_air} use a planar Neumann
penalty with unit weight. Diagnostic residuals use the stored
coefficients in float64.

Reference grids use $192$, $288$, and $384$ cells in each direction with
angular wrapping; the fine array has $385\times385\times384$ values.
The reference stencil requires no diffusion beyond $\nu$ at these
positive-viscosity settings. Grid/time refinement and changes in viscosity or terminal
smoothing are separated in Table~\ref{tab:air_reference}.
The smoothed terminal function and homogeneous Neumann data are
incompatible at $t=T$; Table~\ref{tab:air_boundary} reports the boundary
and residual errors. The neural stencil does not satisfy
\eqref{eq:main_monotonicity} on the full Air3D domain.
\begin{table}[!htbp]
\centering\footnotesize\resizebox{\textwidth}{!}{\begin{tabular}{lccc}
\toprule
Method & Domain residual RMS & FD--AD discrepancy & Boundary normal RMS \\
\midrule
Direct PINN & $0.1331\pm0.0086$ & $0.0211\pm0.0023$ & $0.1224\pm0.0044$ \\
FD-direct & $0.1397\pm0.0055$ & $0.0201\pm0.0015$ & $0.1263\pm0.0085$ \\
AD-PINN-PI & $0.1454\pm0.0050$ & $0.0180\pm0.0012$ & $0.1322\pm0.0054$ \\
FD-PINN-PI & $0.1501\pm0.0041$ & $0.0176\pm0.0009$ & $0.1321\pm0.0042$ \\
\bottomrule
\end{tabular}}
\caption{Residual and boundary errors at the fixed update
budget. Entries are sampled RMS means and sample SDs over five paired seeds.}
\label{tab:air_boundary}
\end{table}
\begin{table}[!htbp]
\centering\footnotesize\resizebox{\textwidth}{!}{\begin{tabular}{llrr}
\toprule
Channel & Variation & Relative difference [\%] & Max sampled difference \\
\midrule
Numerical refinement & 192 vs 384 & 0.01549 & 0.000457 \\
Numerical refinement & 288 vs 384 & 0.00446 & 0.000182 \\
Numerical refinement & Half time step & 0.00268 & 0.000077 \\
Terminal smoothing & $\epsilon=0$ & 0.72032 & 0.016753 \\
Terminal smoothing & $\epsilon=0.05$ & 0.53441 & 0.012235 \\
Terminal smoothing & $\epsilon=0.2$ & 1.99384 & 0.041570 \\
Artificial viscosity & $\nu=0.1$ & 3.97289 & 0.110494 \\
Artificial viscosity & $\nu=0.4$ & 7.44777 & 0.160325 \\
\bottomrule
\end{tabular}}
\caption{Air3D reference sensitivities at $t=0$, relative to the
$384$-cell, $\nu=0.2$, $\epsilon=0.1$ baseline. The same independent
central-box states give relative differences and sampled maxima.}
\label{tab:air_reference}
\end{table}

The FD-increment comparison uses three paired seeds at each increment,
keeping the prescribed update budget and target fixed. In
Table~\ref{tab:air_increments}, consistency is the RMS difference between
FD and continuum residuals, evaluated in float64; precision is the RMS
difference between float32 and float64 FD evaluation of the same network.
\begin{table}[!htbp]
\centering\footnotesize\setlength{\tabcolsep}{3pt}\begin{tabular}{rrlrrr}
\toprule
$h$ & $\tau$ & Method & Error (\%) & \shortstack{Consistency\\($10^{-3}$)} & \shortstack{Precision\\($10^{-3}$)} \\
\midrule
0.025 & 0.000625 & FD-direct & $1.275 \pm 0.127$ & $5.290 \pm 0.548$ & $1.694 \pm 0.154$ \\
0.025 & 0.000625 & FD-PINN-PI & $1.421 \pm 0.104$ & $4.643 \pm 0.353$ & $1.891 \pm 0.169$ \\
0.050 & 0.002500 & FD-direct & $1.284 \pm 0.047$ & $20.224 \pm 2.144$ & $0.446 \pm 0.011$ \\
0.050 & 0.002500 & FD-PINN-PI & $1.442 \pm 0.046$ & $17.026 \pm 0.465$ & $0.474 \pm 0.039$ \\
0.100 & 0.010000 & FD-direct & $1.289 \pm 0.230$ & $78.996 \pm 9.832$ & $0.098 \pm 0.015$ \\
0.100 & 0.010000 & FD-PINN-PI & $1.386 \pm 0.086$ & $71.938 \pm 6.881$ & $0.121 \pm 0.008$ \\
\bottomrule
\end{tabular}
\caption{Air3D increment sensitivity: means and sample SDs over paired
seeds $0,1,2$. Errors are relative percentages; consistency and precision
are residual RMS differences, reported in units of $10^{-3}$.}
\label{tab:air_increments}
\end{table}

The reference with $\nu=0$ includes the local-Lax--Friedrichs
stabilization required on its grid. The comparison in
Table~\ref{tab:air_zero} uses the same terminal smoothing and common
central-box states. Enlarging the planar half extent from $3$ to $4.5$
and $6$ keeps spatial spacing and the time step fixed, with $384$ angular
cells, $\nu=0.2$, and no added diffusion. Relative to half extent $6$,
the pooled difference is $0.00187\%$ at half extent $3$ and is at
floating-point precision at half extent $4.5$.
\begin{table}[!htbp]
\centering\small\begin{tabular}{rrcc}
\toprule
$\nu$ & Cells & $t=0$ difference [\%] & Pooled difference [\%] \\
\midrule
0.05 & 720 & 1.879 & 1.325 \\
0.1 & 384 & 4.013 & 2.820 \\
0.2 & 384 & 8.077 & 5.669 \\
0.4 & 384 & 15.686 & 11.019 \\
\bottomrule
\end{tabular}
\caption{Positive-viscosity references versus the $720$-cell reference
with $\nu=0$. Terminal smoothing remains $0.1$.}
\label{tab:air_zero}
\end{table}

\section{Prescribed errors and numerical refinement}\label{subsec:periodic_defects}
For the exact solution and coefficients in Section~\ref{sec:main_defects}, set
$\widetilde v^0=v_\star+(T-t)\sin(x_1-x_2)$. At iteration $n$, the
momentum perturbation is
$m_n(\cos(2x_1+x_2),-\sin(x_1-x_2))$.
The minimizing feedback becomes
$-0.12(1-2a_n)\operatorname{sign}(p_1)$, with exact control optimization error
$0.24a_n|p_1|$. Terminal data are changed by
$z_n\sin(2x_1+x_2)$. The amplitudes are
$m_n=0.12s_n$, $a_n=0.15s_n$, $z_n=0.04s_n$, where $s_n$ equals
$1$, $0.55^{n-1}$, or $\mathbf1_{\{n=2\}}$; inactive perturbations are zero.

The reconstruction is Fourier in space and piecewise linear in time.
Feedbacks are recomputed from the reconstructed preceding momentum on
shifted validation grids. The continuum policy-evaluation residual is
computed from the reconstruction and includes discretization error.
The bounds $B_n$ and $B_n^r$ in~\eqref{eq:main_diagnostic_history} use
unweighted defect integrals and the weighted initial gradient energy.
In the numerical evaluation, spatial integrals are divided by $|\Omega|$,
accounting for the single volume factor in that formula.

The common refinement matrix uses $33^2$, $65^2$, and $129^2$ spatial
nodes with $128$, $256$, and $512$ RK4 time intervals. From $65^2$ to
$129^2$, the final persistent momentum energy norm changes by $1.67\%$;
the other three cases change by less than $0.0014\%$.
Independent Fourier identities and exact quadratic-time integration
reproduce the value, momentum, terminal, and control optimization norms in floating
point. The sign-dependent evaluation residual still uses quadrature.

At $1024$ time intervals, increasing from $129^2$ to $257^2$ nodes
changes the final persistent-momentum energy norm by $0.167\%$;
time and quadrature refinement have smaller effects.
Section~\ref{sec:main_defects} reports all twelve cases from the common matrix.

\section*{Acknowledgments}
Yeoneung Kim is supported by the National Research Foundation of Korea (NRF)
grants funded by the Korea government (MSIT) (RS-2023-00219980,
RS-2026-25490291) and Yonsei University.

Under the authors' direction, OpenAI Codex assisted with manuscript
organization and editing, code development, and mathematical and numerical
consistency checks. The authors take full responsibility for the content
of this work.

\section*{Data and code availability}
The experimental code, configurations, and table and figure generators are available at
\url{https://github.com/yeoneung/fdm_pinn}.
The repository README documents the \texttt{run.py run-all},
\texttt{run.py analyze}, and \texttt{run.py plot} commands, from experiment
execution to checkpoint evaluation and generation of the tables and figures.
The analysis and plotting routines are in \texttt{analysis/} and
\texttt{plots/}, respectively.

\bibliographystyle{unsrtnat}
\bibliography{refs}
\end{document}